\documentclass[10pt,twocolumn,fleqn]{article}

\usepackage{amssymb}
\usepackage{ifthen}
\usepackage{amsmath}
\usepackage{stmaryrd}
\usepackage{bm}
\usepackage{pgfplots}

\newcommand{\R}{\ensuremath{\mathbb{R}}}
\newcommand{\mat}[1]{\ensuremath{\begin{pmatrix}#1\end{pmatrix}}}
\newcommand{\RM}[1]{\ensuremath{\MakeUppercase{\romannumeral #1}}}
\newcommand{\Div}[2][]{\ensuremath{\,\text{div}_{#1}(#2)}}
\newcommand{\cof}[1]{\ensuremath{\,\text{cof}(#1)}}
\newcommand{\T}{\ensuremath{\mathcal{T}}}
\newcommand{\E}{\ensuremath{\mathcal{E}}}
\newcommand{\Jbnd}{\ensuremath{J_{b}}}
\newcommand{\tr}[1]{\ensuremath{\,\text{tr}(#1)}}
\newcommand{\Av}[1]{\ensuremath{\{#1\}}}
\newcommand{\re}[1]{\ensuremath{\hat{#1}}}
\newcommand{\ph}[1]{\ensuremath{#1}}
\newcommand{\HDivDiv}[1][]{\ensuremath{H(\text{divdiv}\ifthenelse{\equal{#1}{}}{}{,#1})}}
\newcommand{\HDiv}[1][]{\ensuremath{H(\text{div}\ifthenelse{\equal{#1}{}}{}{,#1})}}
\newcommand{\Hone}[1][]{\ensuremath{H^1\ifthenelse{\equal{#1}{}}{}{(#1)}}}
\newcommand{\Cone}[1][]{\ensuremath{C^1\ifthenelse{\equal{#1}{}}{}{(#1)}}}
\newcommand{\Ltwo}[1][]{\ensuremath{L^2\ifthenelse{\equal{#1}{}}{}{(#1)}}}
\newcommand{\idop}{\text{id}}
\newcommand{\bsigma}{\bm{\sigma}}

\newcommand{\refshell}{\re{\mathcal{S}}}
\newcommand{\physshell}{\ph{\mathcal{S}}}
\newcommand{\refT}{\re{\mathcal{T}}}
\newcommand{\physT}{\ph{\mathcal{T}}}
\newcommand{\refE}{\re{\mathcal{E}}}
\newcommand{\physE}{\ph{\mathcal{E}}}
\newcommand{\physmom}{\ph{\bsigma}}

\newcommand{\refnv}{\re{\nu}}
\newcommand{\physnv}{\ph{\nu}}
\newcommand{\refbnv}{\re{\mu}}
\newcommand{\physbnv}{\ph{\mu}}
\newcommand{\reft}{\re{\tau}_e}
\newcommand{\physt}{\ph{\tau}_e}
\newcommand{\reftang}{\re{\tau}}
\newcommand{\phystang}{\ph{\tau}}

\newcommand{\Dim}{3}
\newcommand{\Lagr}{\mathcal{L}}
\newcommand{\physgrad}{\nabla_{\phystang}}
\newcommand{\refgrad}{\nabla_{\reftang}}
\newcommand{\defgrad}{\bm{F}}
\newcommand{\defhesse}{\bm{H}}
\newcommand{\Imat}{\bm{I}}
\newcommand{\Mmat}{\bm{M}}
\newcommand{\Amat}{\bm{A}}
\newcommand{\Ptangref}{\bm{P}_{\reftang}}

\newcommand{\Ptauphys}{\bm{P}_{\physt}}
\newcommand{\CGstrain}{\bm{C}}
\newcommand{\Gstrain}{\bm{E}}
\newcommand{\Hessian}{\bm{\mathcal{H}}}

\numberwithin{equation}{section}
\usepackage{chngcntr}
\counterwithin{figure}{section}
\counterwithin{table}{section}

\begin{document}
\begin{center}
{\Large \textbf{The Hellan--Herrmann--Johnson Method for Nonlinear Shells}} \\[0.5cm]
\end{center}
Michael Neunteufel$^{a,}$\footnote{Corresponding author.\\E-mail adress: michael.neunteufel@tuwien.ac.at }, Joachim Sch{\"o}berl$^{a}$\\
$^a$ Institute for Analysis and Scientific Computing, TU Wien, Wiedner Hauptstra{\ss}e 8-10, 1040 Wien, Austria

\begin{abstract}
In this paper we derive a new finite element method for nonlinear shells. The Hellan--Herrmann--Johnson (HHJ) method is a mixed finite element method for fourth order Kirchhoff plates. It uses convenient Lagrangian finite elements for the vertical deflection, and introduces sophisticated finite elements for the moment tensor. In this work we present a generalization of this method to nonlinear shells, where we allow finite strains and large rotations. The geometric interpretation of degrees of freedom allows a straight forward discretization of structures with kinks. The performance of the proposed elements is demonstrated by means of several established benchmark examples.\newline

\textbf{\textit{Keywords:}} nonlinear shells; structural mechanics; discrete differential geometry; mixed finite elements; Kirchhoff hypothesis
\end{abstract}
 
\section{Introduction}
\label{sec:introduction}
The difficulty of constructing simple $\Cone$-conforming Kirchhoff--Love shell elements led to the development of the well-known discrete Kirchhoff (DKT) elements \cite{Mor71,KB96,BZH01}, where the Kirchhoff constraint was inforced in a discrete way along the edges. The class of rotation-free (RF) elements eliminate the rotational degrees of freedom by using out-of-plane translation degrees of freedom (dofs) \cite{OZ00,BS06,GT07}. Alternative approaches are discontinuous Galerkin (DG) methods \cite{EGHL02,HLM2017,VPS2017} and Isogeometric Analysis (IGA) \cite{HCB05,SF18,KBLW09,EOB2013}.

The HHJ method for fourth order Kirchhoff plates has been developed and analyzed in \cite{Hel67,Her67,Joh73}. Later work has been done in the 80s \cite{Com89,AB85}, 90s \cite{St91} and recently after 20 years \cite{Chen18,Huang11,BPS17}. It overcomes the issue of $\Cone$-conformity by introducing the moment tensor as an additional tensor field leading to a mixed method. The tangential displacement and normal-normal stress method (TDNNS) developed for linear elasticity and Reissner--Mindlin plates in \cite{PS11,PS12,PS17,PS18} follows the idea of mixed methods, where the stress tensor gets interpolated in the reinvented $\HDivDiv$ space from the HHJ method.

In this paper modern coordinate-free differential geometry, see e.g. \cite{DZ11,Sp79}, is used to define the shell energy. The aim of this work is to find a (high-order) finite element shell element, consisting of $\Hone$-conforming finite elements for the displacement and $\HDivDiv$ elements for the moments. It turns out, that this model can be seen as a generalization of the HHJ method to nonlinear shells. Furthermore, the method can handle surfaces with kinks in a natural way without additional treatment. Numerical results are shown to confirm the model.

\section{Methodology}
\label{sec:methodology}
\subsection{Notation and finite element spaces}
\label{subsec:notation_fe_spaces}
Let $\physshell$ be a 2-dimensional surface in $\R^3$, and let $\physshell_h=\bigcup_{T\in\T_h}T$ be its approximation by a triangulation $\physT_h$ consisting of possibly curved triangles or quadrilaterals. The set of all edges in $\physT_h$ is denoted by $\physE_h$. Further, let $\Ltwo[\physshell_h]$ and $C^0(\physshell_h)$ be the set of all square-integrable and continuous functions on $\physshell_h$, respectively. 

\begin{figure}[h!]
\centering
\includegraphics[width=0.4\textwidth]{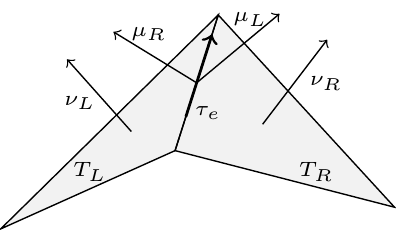}
\caption{Normal, element-normal and normalized edge tangent vectors on two triangles $T_L$ and $T_R$.}
\label{fig:v_conv_tang_trig}
\end{figure}

For each element in $\physT_h$ we denote the surface normal vector by $\physnv$ and the  normalized edge tangent vector between elements by $\physt$. The outgoing element-normal vector $\physbnv$ is defined as $\physbnv=\pm\physnv\times\physt$ depending on the orientation of $\physt$, cf. Figure \ref{fig:v_conv_tang_trig}.

The set of all piece-wise polynomials of degree $k$ on $\T_h$ is denoted by $\Pi^k(\physT_h)$. With this, we define the following function and finite element spaces
\begin{flalign}
&\Hone[\T_h]:=\left\{u \in \Ltwo[\physshell_h]\,|\, \physgrad u\in[\Ltwo[\physshell_h]]^3\right\},\label{eq:h1_fs}&&\\
&\Sigma(\T_h):=\left\{\physmom \in [C^{\infty}(\T_h)]^{\Dim\times \Dim}_{sym}\,|\, \llbracket\physmom_{\physbnv\physbnv}\rrbracket=0\right\},\label{eq:hdivdiv_fs}&&\\
&V^k_h(\T_h):=\Pi^k(\T_h)\cap C^0(\physshell_h),\label{eq:h1_fes}&&\\
&\Sigma_h^k(\T_h):=\left\{\physmom \in [\Pi^k(\T_h)]^{\Dim\times \Dim}_{sym}\,|\, \llbracket\physmom_{\physbnv\physbnv}\rrbracket=0\right\},\label{eq:hdivdiv_fes}&&\\
&\overline{\Gamma}_h^k(\T_h):=\left\{u \in [\Pi^k(\T_h)]^{\Dim}\,|\, \llbracket u_{\physbnv} \rrbracket=0\right\},\label{eq:hdiv_fes}&&
\end{flalign}
where we used the notations $\physmom_{\physbnv\physbnv}:=\physbnv^T\physmom \physbnv$ and $u_\physbnv :=u\cdot \physbnv$ with $\llbracket\cdot\rrbracket$ denoting the jump over elements. Note, that $\physgrad u$ denotes the surface gradient of $u$, which can be introduced in weak sense \cite{de13}, or directly as Fr\'{e}chet-derivative.

For the construction of finite element spaces and explicit basis functions of \eqref{eq:h1_fes}, \eqref{eq:hdivdiv_fes} and  \eqref{eq:hdiv_fes} we refer to \cite{Br2013,PS11,Zaglmayr06,Nedelec1986,Brezzi1985,ZT20}. With a hierarchical basis of \eqref{eq:hdiv_fes} we can define the finite element space $\Gamma_h^k(\T_h)$ as the space $\overline{\Gamma}_h^k(\T_h)$, where the inner degrees of freedom are neglected. $\Sigma_h^k(\T_h)$ will be called the Hellan--Herrmann--Johnson finite element space.\newline

In the following we denote the Frobenius scalar product of two matrices $\bm{A}$, $\bm{B}$ by $\bm{A}:\bm{B}:=\sum_{i,j}\bm{A}_{ij}\bm{B}_{ij}$, $\|\bm{A}\|_F:=\sqrt{\bm{A}:\bm{A}}$ and $\sphericalangle(a,b):=\arccos(\frac{a\cdot b}{\|a\|_2\|b\|_2})$ measures the angle between two vectors $a$, $b$, $\|\cdot\|_2$ denoting the Euclidean norm.

\subsection{Shell model}
\label{subsec:shell_model}
Let $\re{\Omega}\subset\R^3$ be an undeformed configuration of a shell with thickness $t$, described by the mid-surface $\refshell$ and the according orientated normal vector $\refnv_{\refshell}$
\begin{flalign}
\re{\Omega}:=\{\re{x}+z\refnv_{\refshell}(\re{x}): \re{x}\in \refshell, z\in [-t/2,t/2]\}.&&
\end{flalign}
Furthermore, let $\Phi:\re{\Omega}\rightarrow\ph{\Omega}$ be the deformation from the initial to the deformed configuration of the shell and $\phi:\refshell_h\rightarrow\physshell_h$ the deformation of the approximated mid-surface. I.e., let $\phi\in [V_h^{k+1}(\refT_h)]^3$ with $\refT_h$ and $\physT_h=\phi(\refT_h)$ the according triangulations of $\refshell_h$ and $\physshell_h$. Then, we define $\defgrad:=\refgrad\phi$ and $J:=\|\cof{\defgrad}\|_F=\|\cof{\defgrad}\refnv\|_2$ as the deformation gradient and the deformation determinant, respectively. Here, $\cof{\defgrad}$ denotes the cofactor matrix of $\defgrad$. We can split the deformation into the identity function and the displacement, $\phi=\idop+u$, and thus, $\defgrad= \Ptangref+\refgrad u$ with the projection onto the tangent plane $\Ptangref:=\Imat-\refnv\otimes\refnv$, $\otimes$ denoting the dyadic outer product.

We consider the Kirchhoff--Love assumption, where the deformed normal vector has to be orthogonal to the deformed mid-surface $\physshell_h$. With Steiner's formula, asymptotic analysis in the thickness parameter $t$ and using the plane strain assumption for the material norm, we obtain for the according shell energy functional 
\begin{flalign}
\mathcal{W} = \frac{t}{8}\|\ph{\RM{1}}-\re{\RM{1}}\|_{\Mmat}^2+\frac{t^3}{24}\|\ph{\RM{2}}-\re{\RM{2}}\|^2_{\Mmat}.&&\label{eq:shell_energy_start}
\end{flalign}
$\mathcal{W}$ is given in terms of differential forms, see \cite{SIMO89}, and \eqref{eq:shell_energy_start} is comparable to the classical formulations \cite{C05, BRI17, Chap11}.

The material norm is given by
\begin{flalign}
\|\cdot\|^2_{\Mmat}:=\frac{E}{1-\nu^2}\int_{\refshell_h}\left(\nu\tr{\cdot}^2+(1-\nu)\tr{\cdot^2}\right)dx,&&
\end{flalign}
with $E$ the Young's modulus and $\nu$ the Poisson's ratio. $\re{\RM{1}}$, $\ph{\RM{1}}$ and $\re{\RM{2}}$, $\ph{\RM{2}}$ denote the (pull-backed) first and second fundamental form of the reference and deformed configuration, respectively. With the Green strain tensor $\Gstrain:= 1/2(\CGstrain-\Ptangref)$ restricted on the tangent space, $\CGstrain=\defgrad^T\defgrad$ denoting the Cauchy-Green tensor, we obtain
\begin{flalign}
E_{\text{mem}}:=\frac{t}{8}\|\ph{\RM{1}}-\re{\RM{1}}\|_{\Mmat}^2 = \frac{t}{2}\|\Gstrain\|_{\Mmat}^2.\label{eq:shell_energy_mem}&&
\end{flalign}
This corresponds to the membrane energy of the shell. The difference between the curvature of the deformed and initial second fundamental form describes the bending energy, for which holds
\begin{flalign}
\tilde{E}_{\text{bend}}&:= \frac{t^3}{24}\|\ph{\RM{2}}-\re{\RM{2}}\|_{\Mmat}^2&&\nonumber\\ &=\frac{t^3}{24}\|\defgrad^T\refgrad(\physnv\circ\phi)-\refgrad\refnv\|_{\Mmat}^2.&&\label{eq:shell_energy_bend_cont}
\end{flalign}
Motivated by discrete differential geometry, see \cite{GGRZ06} and references therein, and DG methods \cite{DBCM01} we add also distributional contributions to the bending energy
\begin{flalign}
E_{\text{bend}}:= &\frac{t^3}{24}\big(\sum_{\re{T}\in\refT_h}\|\defgrad^T\refgrad(\physnv\circ\phi)-\refgrad\refnv\|_{\Mmat,\re{T}}^2&&\nonumber\\
&+\sum_{\re{E}\in\refE_h}\|\sphericalangle(\physnv_L, \physnv_R)\circ\phi-\sphericalangle(\refnv_L, \refnv_R)\|_{\Mmat,\re{E}}^2\big).&&\label{eq:shell_energy_bend}
\end{flalign}

Thus, with the notation of \eqref{eq:shell_energy_mem} and \eqref{eq:shell_energy_bend}, we have to minimize
\begin{flalign}
\label{eq:nonl_shell_method_first_u}
\tilde{\mathcal{W}}(u) := E_{\text{mem}}+E_{\text{bend}}.&&
\end{flalign}

To reduce this fourth order problem to a second order one, we introduce a new variable $\physmom$ which leads to a mixed saddle point problem. Hence, we have to find the critical points of the following Lagrange functional, which is equivalent to minimize \eqref{eq:nonl_shell_method_first_u}, see Appendix A,
\begin{flalign}
\label{eq:nonl_shell_method_first}
\tilde{\Lagr}(u,\physmom) := \frac{t}{2}\|\Gstrain\|_{\Mmat}^2 -\frac{6}{t^3}\|\physmom\|_{\Mmat^{-1}}^2+\tilde{B}(\physmom,u),&&
\end{flalign}
where
\begin{flalign}
\tilde{B}(\physmom,u):=&\sum_{\re{T}\in\refT_h}\langle\physmom,\defgrad^T\refgrad(\physnv\circ\phi)-\refgrad\refnv\rangle_{\re{T}} &&\nonumber\\
&- \sum_{\re{E}\in\refE_h}\langle\sphericalangle(\physnv_L, \physnv_R)\circ\phi-\sphericalangle(\refnv_L, \refnv_R),\physmom_{\refbnv\refbnv}\rangle_{\re{E}},&&
\end{flalign}
with $\langle\cdot,\cdot\rangle$ denoting the $\Ltwo$-scalar product on an element $\re{T}$ or on an edge $\re{E}$.
With some computations, see Appendix A, we finally obtain the following Lagrange functional
\begin{flalign}
\label{eq:nonl_shell_method}
\Lagr(u,\physmom) = \frac{t}{2}\|\Gstrain\|_{\Mmat}^2 -\frac{6}{t^3}\|\physmom\|^2_{\Mmat^{-1}}-B(\physmom,u),&&
\end{flalign}
with
\begin{flalign}
&B(\physmom,u) =&&\nonumber\\ &\sum_{\re{T}\in\refT_h}\int_{\re{T}}\physmom:(\Hessian_{\physnv\circ\phi} + (1-\refnv\cdot\physnv\circ\phi)\refgrad\refnv)\,dx&&\nonumber\\
&-\sum_{\re{E}\in\refE_h}\int_{\re{E}}(\sphericalangle(\refnv_L, \refnv_R)-\sphericalangle(\physnv_L, \physnv_R)\circ\phi) \physmom_{\refbnv\refbnv}\,ds.\label{eq:final_B}&&
\end{flalign}

In \eqref{eq:final_B} $\Hessian_{\physnv\circ\phi}:=\sum_{i}(\refgrad^2u_i)\physnv_i\circ\phi$, where $\refgrad^2$ denotes the surface Hessian \cite{de13}. For the deformed normal and tangent vectors the following identities hold
\begin{flalign}
&\physnv\circ\phi=\frac{1}{\|\cof{\defgrad}\refnv\|_2}\cof{\defgrad}\refnv=\frac{1}{J}\cof{\defgrad}\refnv,&&\\
&\physt\circ\phi=\frac{1}{\|\defgrad\reft\|_2}\defgrad\reft=\frac{1}{\Jbnd}\defgrad\reft,&&\\
&\physbnv\circ\phi=\pm(\physnv\circ\phi)\times (\physt\circ\phi)= \frac{1}{\|(\defgrad^{\dagger})^T\refbnv\|_2}(\defgrad^{\dagger})^T\refbnv,&&
\end{flalign} 
where $\defgrad^{\dagger}$ denotes the Moore--Penrose pseudo-inverse of $\defgrad$.

The Lagrange multiplier $\physmom$ has the physical meaning of the moment. Note, that the thickness parameter $t$ appears now also in the denominator and the inverse material tensor
\begin{flalign}
\|\cdot\|_{\Mmat^{-1}}^2:= \frac{1+\nu}{E}\int_{\refshell_h}(\tr{\cdot^2}-\frac{\nu}{2\nu+1}\tr{\cdot}^2)\,dx,&&
\end{flalign}
is used.
 
In case of a flat plane \eqref{eq:final_B} becomes
\begin{flalign}
B(\physmom,u) =& \sum_{\re{T}\in\refT_h}\int_{\re{T}}\physmom:\Hessian_{\physnv\circ\phi} \,dx&&\nonumber\\
&-\sum_{\re{E}\in\refE_h}\int_{\re{E}}\sphericalangle(\physnv_L, \physnv_R)\circ\phi \,\physmom_{\refbnv\refbnv}\,ds.\label{eq:final_B_plate}&&
\end{flalign}

A possible simplification of \eqref{eq:final_B} can be achieved by the approximation
\begin{flalign}
\frac{1}{2}\sphericalangle(\physnv_L, \physnv_R) = \Av{\physnv}\cdot\physbnv_L + \mathcal{O}(|\Av{\physnv}\cdot\physbnv_L|^3)\label{eq:approx_angle_comp},&&
\end{flalign}
where $\Av{\physnv}:=\frac{1}{\|\physnv_L+\physnv_R\|_2}(\physnv_L+\physnv_R)$ denotes the averaged normal vector.\newline

The resulting system is a saddle point problem, which would lead to an indefinite matrix after assembling. To overcome this problem, we can use complete discontinuous elements for the moment $\physmom$ and introduce a hybridization variable $\re{\alpha}\in\Gamma_h^k(\refT_h)$ to reinforce the normal-normal continuity of $\physmom$:
\begin{flalign}
\Lagr(u,\physmom,\re{\alpha}) = \frac{t}{2}\|\Gstrain\|_{\Mmat}^2 -\frac{6}{t^3}\|\physmom\|^2_{\Mmat^{-1}}-B(\physmom,u,\re{\alpha}),\label{eq:lagr_alpha}&&
\end{flalign}
where \eqref{eq:final_B} is now given by
\begin{flalign}
&B(\physmom,u,\re{\alpha}) &&\nonumber\\
&=\sum_{\re{T}\in\refT_h}\int_{\re{T}}\physmom:(\Hessian_{\physnv\circ\phi} + (1-\refnv\cdot\physnv\circ\phi)\refgrad\refnv)\,dx&&\nonumber\\
&-\sum_{\re{E}\in\refE_h}\int_{\re{E}}(\sphericalangle(\refnv_L, \refnv_R)-\sphericalangle(\physnv_L, \physnv_R)\circ\phi) \,\langle\langle\physmom_{\refbnv\refbnv}\rangle\rangle\,ds&&\nonumber\\
&\quad+ \int_{\re{E}}\re{\alpha}_{\refbnv} \llbracket\physmom_{\refbnv\refbnv}\rrbracket\,ds,
\end{flalign}
with $\langle\langle\physmom_{\refbnv\refbnv}\rangle\rangle:=1/2(\physmom_{\refbnv_L\refbnv_L}+\physmom_{\refbnv_R\refbnv_R})$. Due to the hybridization variable $\re{\alpha}$, we can use static condensation to eliminate the moment $\physmom$ locally, which leads to a positive definite problem again. The new unknown $\re{\alpha}$ has the physical meaning of the changed angle, the rotation, between two elements. \newline

For the computation of the jump term we use that 
\begin{flalign}
&\sum_{\re{E}\in\refE_h}\int_{\re{E}}\sphericalangle(\refnv_L, \refnv_R)-\sphericalangle(\physnv_L, \physnv_R)\circ\phi\,ds=&&\nonumber\\
&\quad \sum_{\re{T}\in\refT_h}\int_{\partial \re{T}}\sphericalangle(\Av{\refnv}, \refbnv)-\sphericalangle(\Av{\physnv}, \physbnv)\circ\phi\,ds.&&
\end{flalign}

\begin{figure}[h!]
\centering
\includegraphics[width=0.2\textwidth]{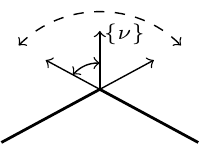}
\includegraphics[width=0.2\textwidth]{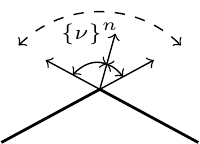}
\caption{Angle computation with the current averaged normal vector $\Av{\physnv}$ and the averaged normal vector $\Av{\physnv}^n$ from the previous step.}
\label{fig:angle_computation}
\end{figure}

To compute the deformed averaged normal vector $\Av{\physnv}$ on an edge, information of the two neighbored elements is needed at once, which would need e.g. Discontinuous Galerkin techniques. Instead, one can use the information of the last (load-step) solution $\Av{\physnv}^{n}$, see Figure \ref{fig:angle_computation}. To measure the correct angle, we have to project $\Av{\physnv}^{n}$ to the plane orthogonal to the tangent vector $\physt$ by using the projection $\Ptauphys^{\perp} = \Imat - \physt\otimes\physt$, and then re-normalize it 
\begin{flalign}
\Av{\physnv}\approx \frac{1}{\|\Ptauphys^{\perp}(\Av{\physnv}^{n})\|_2}\Ptauphys^{\perp}(\Av{\physnv}^{n})=:\overline{\Av{\physnv}}^n.&&\label{eq:approx_averaged_nv}
\end{flalign}

Note that $\physt$ itself depends on the unknown deformation. By using \eqref{eq:approx_averaged_nv} we have to ensure that $\Av{\physnv}^{n}$ lies between the two element-normal vectors, see Figure \ref{fig:angle_computation}. For smooth manifolds the angle between the element-normal vectors tends to $180$ degree as $h\rightarrow0$. Hence, this assumption is fulfilled, if the elements do not rotate more than half of their included angle during one load-step, which is an acceptable and realistic assumption.

\subsection{Relation to the HHJ-method}
\label{subsec:rel_hhj}
If we assume to have a plate which lies in the x-y-plane and a force $f$ is acting orthogonal on it, we can compute the linearized bending energy by solving the following fourth order scalar equation on $\refshell_h$
\begin{flalign}
\Div{\Div{\nabla^2 w}}=f,\label{eq:hhj}&&
\end{flalign}
where the thickness $t$ and all material parameters are hidden in the right-hand side $f$.

Therefore, the HHJ-method \cite{Hel67,Her67,Joh73} introduces the linearized moment tensor $\physmom$ and solves the following saddle point problem instead, given by the Lagrange functional 
\begin{flalign}
\label{eq:hhj_method}
\Lagr(w,\physmom) :=& -\|\physmom\|^2_{\Ltwo[\refshell_h]} + \sum_{\re{T}\in\refT_h}\Big(\int_{\re{T}}\nabla w\cdot\Div{\physmom}\,dx&&\nonumber\\
&-\int_{\partial \re{T}}(\nabla w)_{\physt}\physmom_{\physbnv \physt}\,ds\Big)+\int_{\refshell_h}fw\,dx.&&
\end{flalign}

If we now consider our shell model \eqref{eq:nonl_shell_method}, neglect the membrane energy term and the material parameters and linearize the bending energy, see Appendix B, we obtain \eqref{eq:hhj_method}. Thus, \eqref{eq:nonl_shell_method} can be seen as a generalization of the HHJ-method \eqref{eq:hhj_method} from linear plates to nonlinear shells.

\subsection{Boundary conditions and kink structures}
\label{subsec:bc_kinks}
For $\Hone$ the Dirichlet boundary condition $u=u_D$ can be used to prescribe the displacement on the boundary, whereas the do-nothing condition is used for free boundaries. For $\physmom\in\HDivDiv$ we can prescribe the normal-normal component, $\physmom_{\physbnv\physbnv}$, on the boundary. Homogeneous Dirichlet data, $\physmom_{\physbnv\physbnv}=0$, are used for free boundaries. By setting non-homogeneous data one can prescribe a moment. The do-nothing Neumann boundary condition $\physmom_{\physbnv\physt}=0$ is used for clamped boundaries.

In the case of a complete discontinuous moment tensor and the hybridization variable $\re{\alpha}$, the boundary conditions for $\physmom$ have to be incorporated in terms of $\re{\alpha}$. Note that the essential and natural boundary conditions swap, i.e. the clamped boundary condition is now set directly as homogeneous Dirichlet data and the prescribed moment is handled natural as a right-hand side.

If we compute the variations of \eqref{eq:nonl_shell_method} with respect to $\physmom$, we obtain in strong form that the angle from the initial configuration gets preserved, see \eqref{eq:angle_pres}. The hidden interface condition for the displacement $u$ in strong form are not needed for the method itself. However, if one uses e.g. Residual error estimators, the boundary conditions are crucial, see Appendix C for the calculations.\newline

The method can also handle non-smooth surfaces with kinks and branching shells, where one edge is shared by more than two elements, in a natural way, without any extra treatments. Due to the normal-normal continuity of $\physmom$ the moment gets preserved over the kinks and as the angle is the same on the initial and deformed configuration, the kink itself gets also preserved. Note, that in this case simplification \eqref{eq:approx_angle_comp} cannot be used any more, as $|\Av{\physnv}\cdot\physbnv| \nrightarrow 0$ as $h\rightarrow 0$ at the kinks.

\subsection{Shell element}
\label{subsec_shell_el_assembling}
\begin{figure}[h!]
	\centering
	\includegraphics[width=0.16\textwidth]{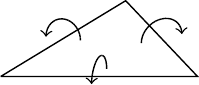}
	\includegraphics[width=0.14\textwidth]{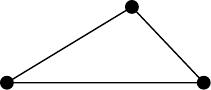}
	\includegraphics[width=0.13\textwidth]{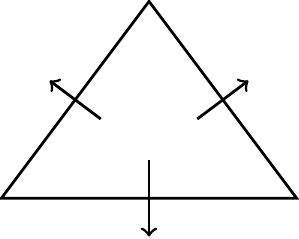}
	
	\includegraphics[width=0.15\textwidth]{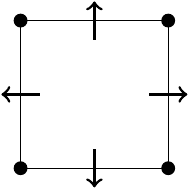}
	\includegraphics[width=0.15\textwidth]{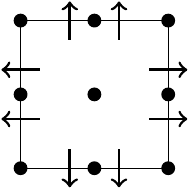}
	\caption{Lowest order $\HDivDiv$, $\Hone$ and $\HDiv$ elements for the moment, displacement and hybridization variable (top) and lowest order and high order hybridized quadrilateral shell element (bottom).}
	\label{fig:fes_shell_el}
\end{figure}

Combining the displacement $u\in [V_h^{k}(\refT_h)]^3$, the moment tensor $\physmom\in\Sigma^{k-1}_h(\refT_h)$ and, eventual, the hybridization space $\Gamma_h^{k-1}(\refT_h)$ leads to our shell element. For polynomial order $k$, the method will be denoted by pk, i.e. p1 is the lowest order method consisting of piece-wise linear displacements and piece-wise constant moments. In Figure \ref{fig:fes_shell_el} the hybridized p1 and p2 element for quadrilaterals can be seen.  Note, that the hybridized lowest order triangle shell element is equivalent to the Morley element \cite{Mor71}. If we use the lowest order elements on triangles for \eqref{eq:nonl_shell_method} then the Hessian term vanishes, as only linear polynomials are used. For quadrilaterals the Hessian is constant on each element in this case.

To solve \eqref{eq:nonl_shell_method} we have to assemble the according matrix. As it is formulated in terms of a Lagrange functional, the first variations must be computed, which is a bit challenging due to the non-linearity but doable, see Appendix C. If, however, the finite element software supports energy based integrators where the variations are calculated automatically, one can use directly the Lagrange functional \eqref{eq:nonl_shell_method}.

\subsection{Membrane locking}
\label{subsec:mem_locking}
We observed that the lowest order elements do not suffer from locking, but for the higher order methods membrane locking, cf. \cite{Pit92}, may occur, e.g. in the benchmark cantilever subjected to end moment, section \ref{subsec:cant_end_mom}. To overcome this problem one can interpolate the membrane stress tensor by a $\Ltwo$-projection into a space of reduced dimension, $\|\mathcal{I}^h_{\Ltwo}\Gstrain\|_{\Mmat}^2$. The projection can be incorporated to \eqref{eq:nonl_shell_method} by introducing an auxiliary variable $\bm{R}$ and adding for the displacement $u\in [V^k_h(\refT_h)]^3$ and $\bm{R}\in[\Pi^{k-1}(\refT_h)]^{3\times 3}_{sym}$   
\begin{flalign}
-\frac{1}{2t}\|\bm{R}\|_{\Mmat^{-1}}^2+\langle \bm{R},\Gstrain\rangle&&\label{eq:l2_proj_mem_locking}
\end{flalign}
to the Lagrange functional. As $\bm{R}$ is discontinuous, we can use static condensation to eliminate it locally. This works well for structured quadrilateral meshes and is similar to reduced integration order methods. For triangles, however, the locking is reduced, but still has an impact to the solution. Here, other interpolation operators and spaces have to be used, which is topic of further research.
 
\section{Numerical results}
\label{sec:num_res}
The method is implemented in the NGS-Py interface, which is based on the finite element library Netgen/NGSolve\footnote{www.ngsolve.org} \cite{Sch97,Sch14}.

We will use the lowest order elements p1and also the p3 method as an high-order example. Table \ref{tab:comp_dofs_shell_el} lists the according number of degrees of freedom for each element.

\begin{table}[h!]
\begin{center}
\begin{tabular}{ccccc}
& p1 T & p1 Q & p3 T & p3 Q \\
\hline
dof/el & 12 & 16 & 36 & 48
\end{tabular}
\caption{Number of degrees of  freedom per (hybridized and condensed) element for triangles (T) and quadrilaterals (Q).}
\label{tab:comp_dofs_shell_el}
\end{center}
\end{table}

\subsection{Cantilever subjected to end shear force}
\label{subsec:cant_shear_force}
\begin{figure}[h!]
\centering
\includegraphics[width=0.45\textwidth]{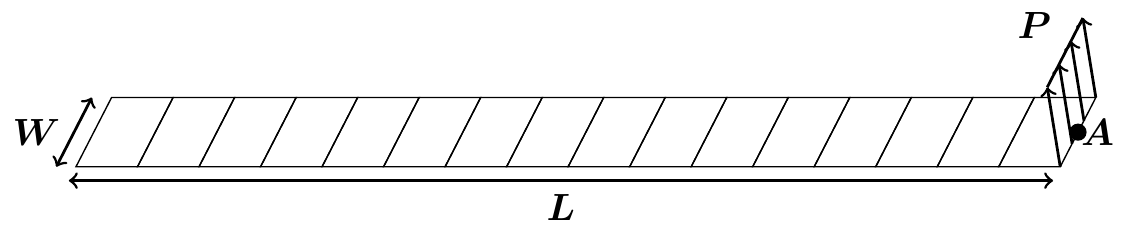}
\caption{Geometry of cantilever subjected to end shear force benchmark.}
\label{fig:res_cant_shear_mesh}
\end{figure}
\begin{figure}[h!]
	\centering
	\includegraphics[width=0.4\textwidth]{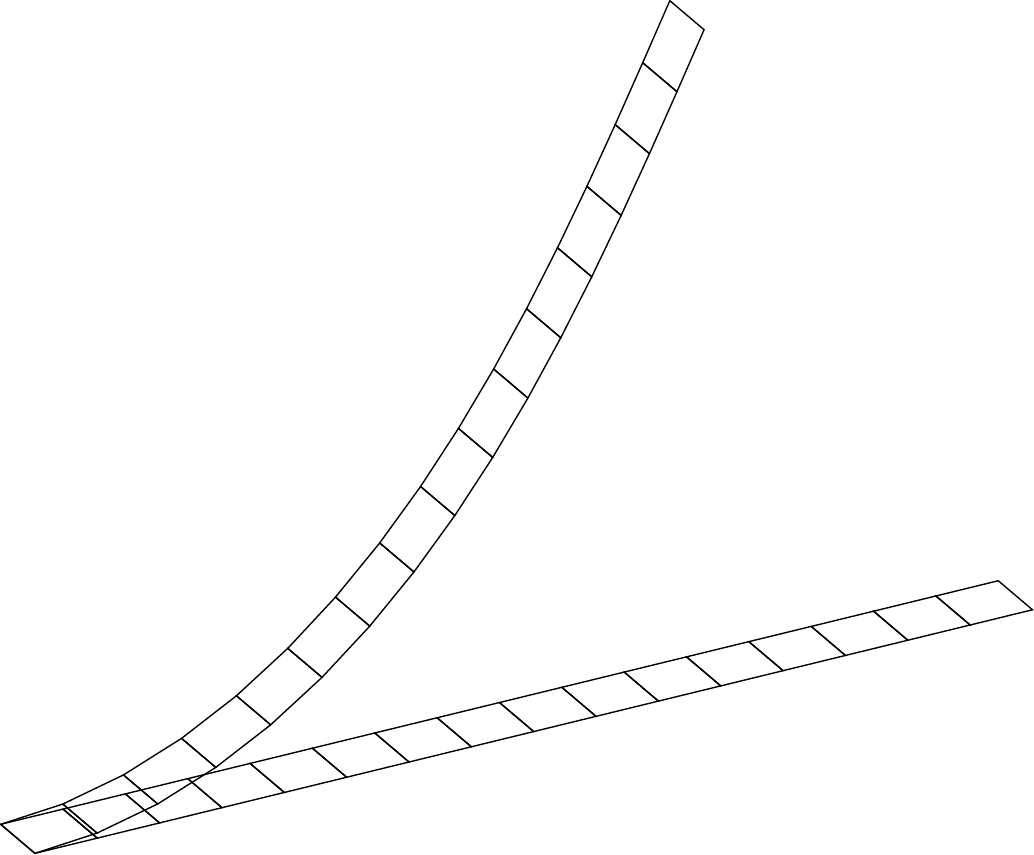}
	\caption{Initial and final configuration of cantilever subjected to end shear force.}
	\label{fig:res_cant_shear_def_mesh}
\end{figure}
\begin{figure}[h!]
\centering
\includegraphics[width=0.45\textwidth]{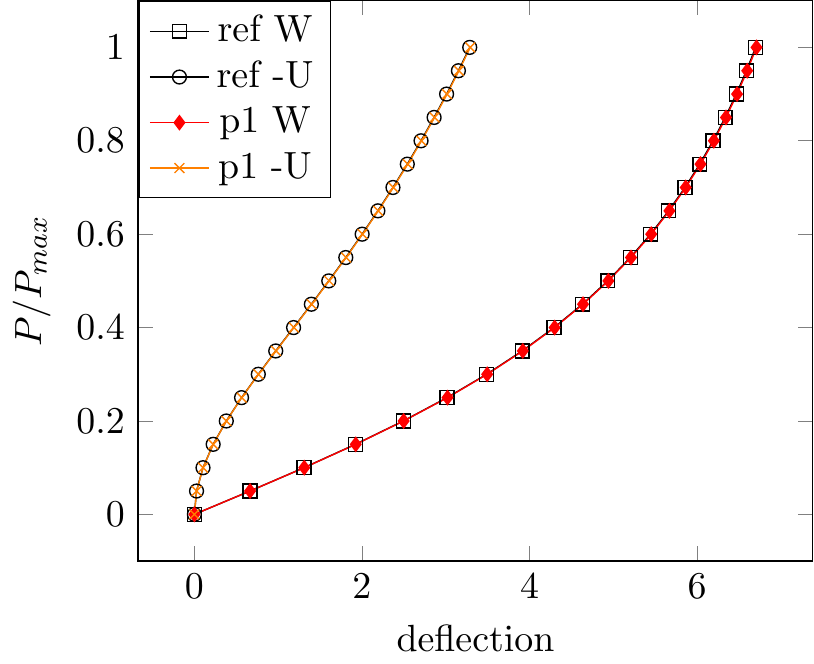}
\caption{Horizontal and vertical load-deflection for cantilever subjected to end shear force with $16\times 1$ grid.}
\label{fig:res_cant_shear}
\end{figure}

An end shear force $P$ on the right boundary is applied to a cantilever, which is fixed on the left. The material and geometrical properties are $E=1.2\times10^6$, $\nu=0$, $L=10$, $W=1$, $t=0.1$ and $P_{max}=4$, see Figure \ref{fig:res_cant_shear_mesh}. A structured $16\times 1$ rectangular grid  is used. The reference values are taken from \cite{Sze04}. In Figure \ref{fig:res_cant_shear_def_mesh} one can see the initial and deformed mesh and in Figure \ref{fig:res_cant_shear} and Table \ref{tab:res_cant_shear} the results.

\begin{table}[h!]
\begin{center}
\begin{tabular}{c|c|c||c|c|c}
	$P/P_{max}$ & -U & V & $P/P_{max}$ & -U & V \\
	\hline
	0.05 & 0.026 & 0.664 & 0.55 & 1.811 & 5.210 \\
	0.10 & 0.104 & 1.311 & 0.60 & 2.007 & 5.452 \\
	0.15 & 0.225 & 1.926 & 0.65 & 2.195 & 5.669 \\
	0.20 & 0.382 & 2.498 & 0.70 & 2.375 & 5.864 \\
	0.25 & 0.565 & 3.021 & 0.75 & 2.546 & 6.040 \\
	0.30 & 0.765 & 3.494 & 0.80 & 2.710 & 6.199 \\
	0.35 & 0.974 & 3.919 & 0.85 & 2.867 & 6.344 \\
	0.40 & 1.187 & 4.299 & 0.90 & 3.015 & 6.476 \\
	0.45 & 1.399 & 4.638 & 0.95 & 3.157 & 6.597 \\
	0.50 & 1.608 & 4.940 & 1.00 & 3.292 & 6.708
\end{tabular}
\caption{Horizontal and vertical deflection of cantilever subjected to end shear force with p1 and $16\times 1$ grid.}
\label{tab:res_cant_shear}
\end{center}
\end{table}

\subsection{Cantilever subjected to end moment}
\label{subsec:cant_end_mom}
\begin{figure}[h!]
	\centering
	\includegraphics[width=0.45\textwidth]{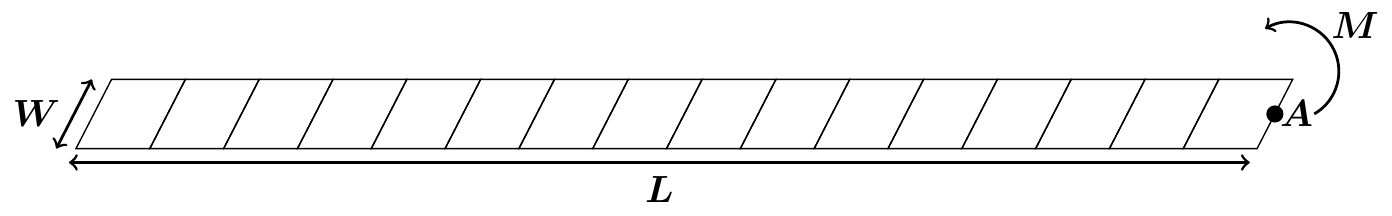}
	\caption{Geometry of cantilever subjected to end moment benchmark.}
	\label{fig:res_cant_bend_mesh}
\end{figure}
\begin{figure}[h!]
	\centering
	\includegraphics[width=0.4\textwidth]{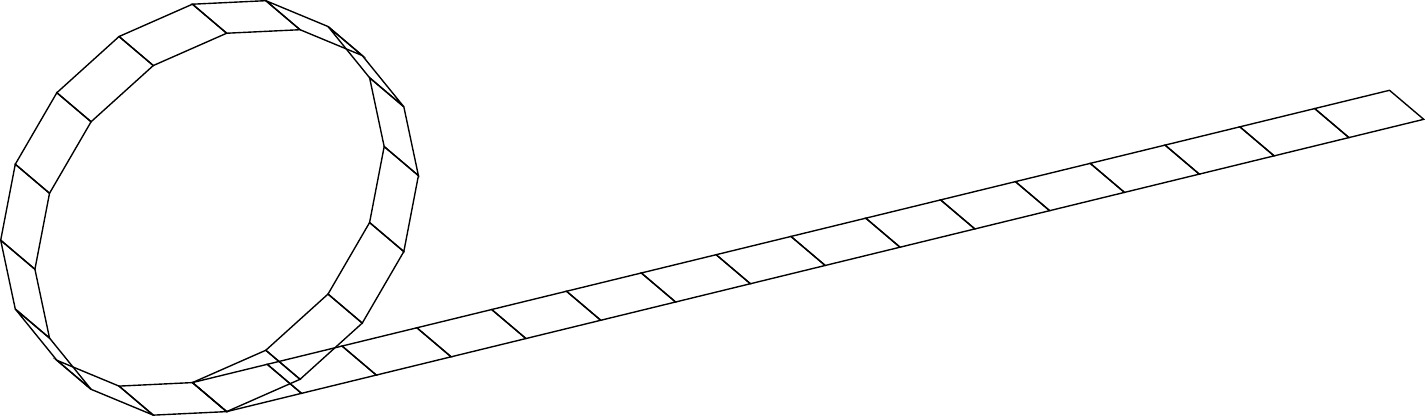}
	\caption{Initial and final configuration of cantilever subjected to end moment.}
	\label{fig:res_cant_bend_def_mesh}
\end{figure}
\begin{figure}[h!]
	\centering
	\includegraphics[width=0.45\textwidth]{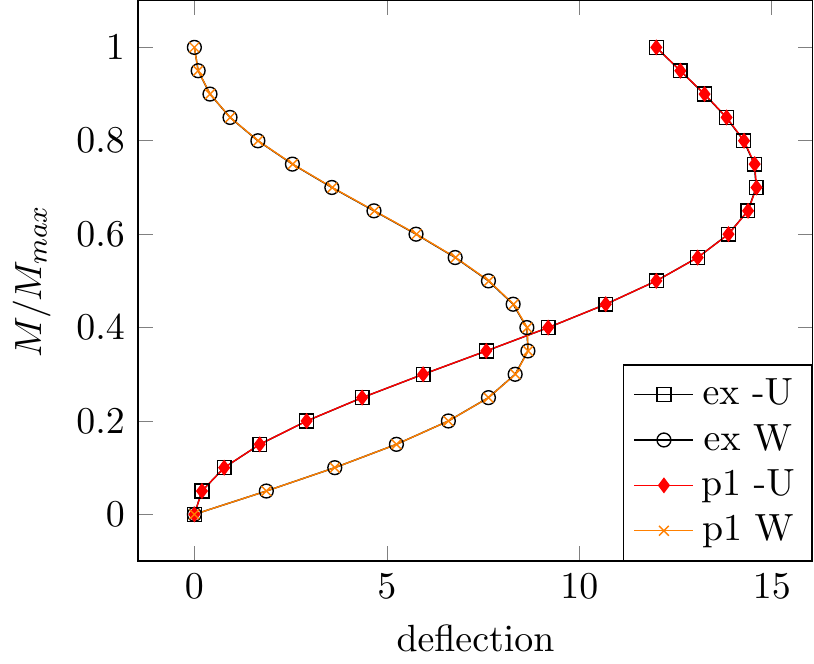}
	\caption{Horizontal and vertical load-deflection for cantilever subjected to end moment with $16\times 1$.}
	\label{fig:res_cant_bend}
\end{figure}

A cantilever is clamped on the left side and a moment $M$ is applied on the right. On the other boundaries we use the symmetry-condition. The material and geometrical properties are $E=1.2\times10^6$, $\nu=0$, $L=12$, $W=1$, $t=0.1$ and $M_{max}=50\pi/3$, see Figure \ref{fig:res_cant_bend_mesh}. The results can be found in Figure \ref{fig:res_cant_bend} and Table \ref{tab:res_cant_bend}, and the initial and final mesh in Figure \ref{fig:res_cant_bend_def_mesh}.

\begin{table}[h!]
\begin{center}
\begin{tabular}{c|c|c|c|c}
	$M/M_{max}$ & $U$ & $U_{ex}$ & W & $W_{ex}$ \\
	\hline
	0.05 & -0.196 & -0.196 & 1.870 & 1.870 \\
	0.10 & -0.773 & -0.774 & 3.648 & 3.648 \\
	0.15 & -1.698 & -1.699 & 5.249 & 5.248 \\
	0.20 & -2.916 & -2.918 & 6.600 & 6.598 \\
	0.25 & -4.357 & -4.361 & 7.643 & 7.639 \\
	0.30 & -5.942 & -5.945 & 8.338 & 8.333 \\
	0.35 & -7.582 & -7.585 & 8.671 & 8.664 \\
	0.40 & -9.191 & -9.194 & 8.646 & 8.637 \\
	0.45 & -10.687 & -10.688 & 8.291 & 8.281 \\
	0.50 & -12.000 & -12.000 & 7.652 & 7.639 \\
	0.55 & -13.075 & -13.073 & 6.788 & 6.775 \\
	0.60 & -13.875 & -13.871 & 5.772 & 5.758 \\
	0.65 & -14.384 & -14.377 & 4.678 & 4.665 \\
	0.70 & -14.603 & -14.595 & 3.583 & 3.571 \\
	0.75 & -14.556 & -14.546 & 2.556 & 2.546 \\
	0.80 & -14.280 & -14.270 & 1.656 & 1.650 \\
	0.85 & -13.826 & -13.818 & 0.931 & 0.926 \\
	0.90 & -13.254 & -13.247 & 0.407 & 0.405 \\
	0.95 & -12.625 & -12.621 & 0.099 & 0.098 \\
	1.00 & -12.000 & -12.000 & 0.000 & 0.000 \\
\end{tabular}
\caption{Horizontal and vertical deflection of cantilever subjected to end moment for p1 and $16\times 1$ grid.}
\label{tab:res_cant_bend}
\end{center}
\end{table}

\subsection{Slit annular plate}
\label{subsec:slit_ann_plate}
\begin{figure}[h!]
	\centering
	\includegraphics[width=0.4\textwidth]{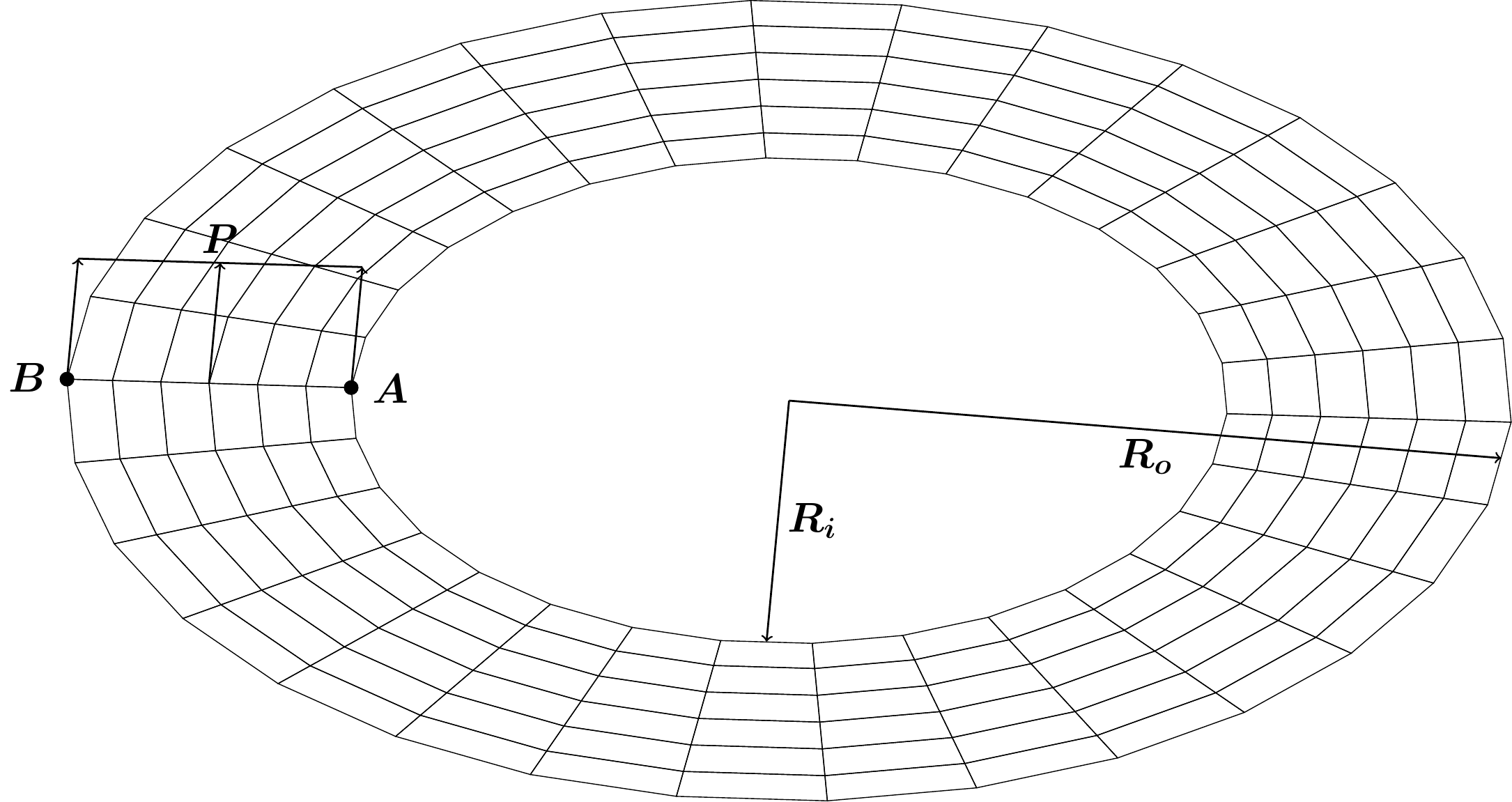}
	\caption{Geometry, force and points of interest of slit annular plate.}
	\label{fig:res_slit_ann_mesh}
\end{figure}

\begin{figure}[h!]
	\centering
	\includegraphics[width=0.4\textwidth]{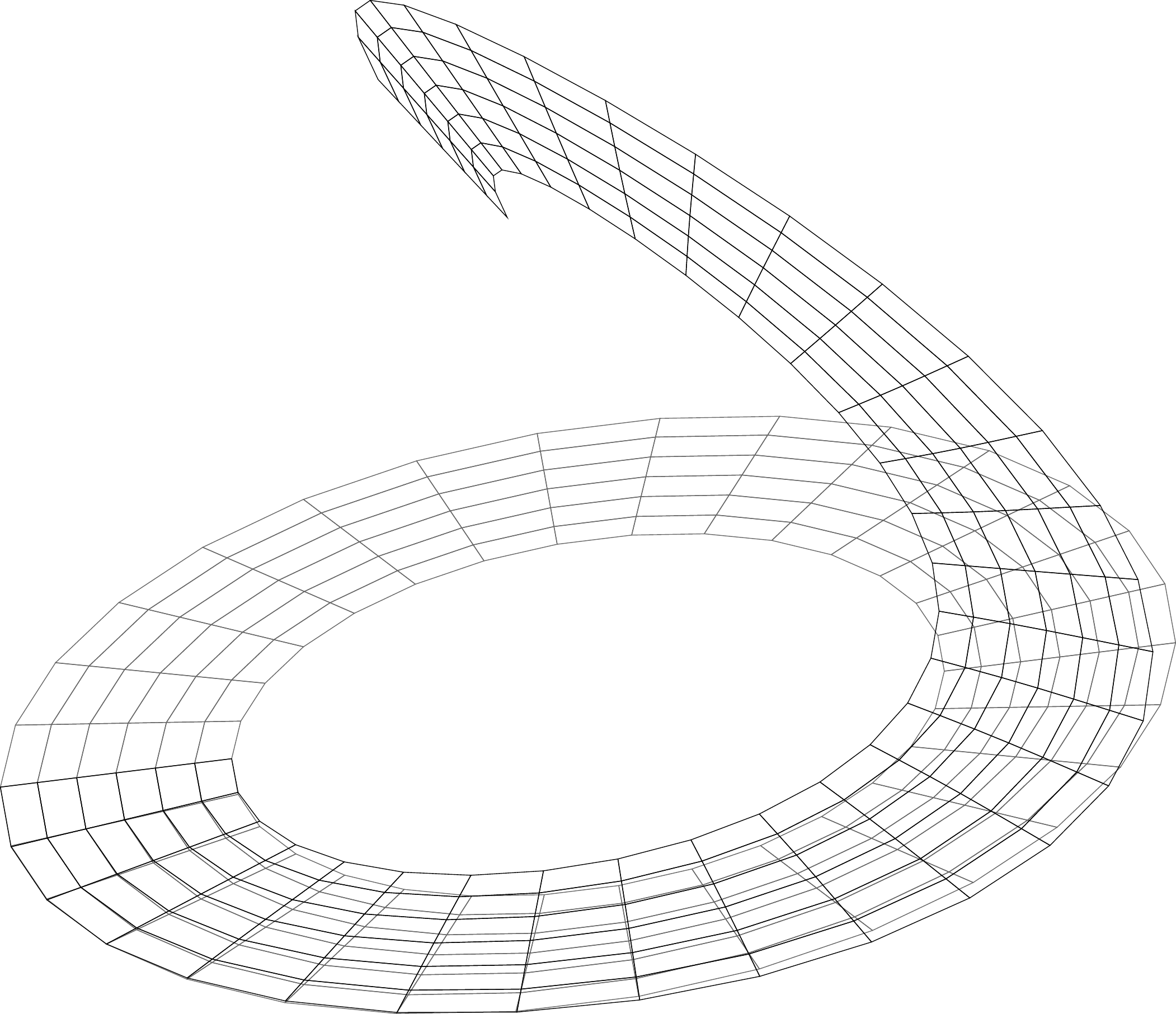}
	\caption{Initial and final configuration of slit annular plate.}
	\label{fig:res_slit_ann_def_mesh}
\end{figure}

\begin{figure}[h!]
	\centering
	\includegraphics[width=0.45\textwidth]{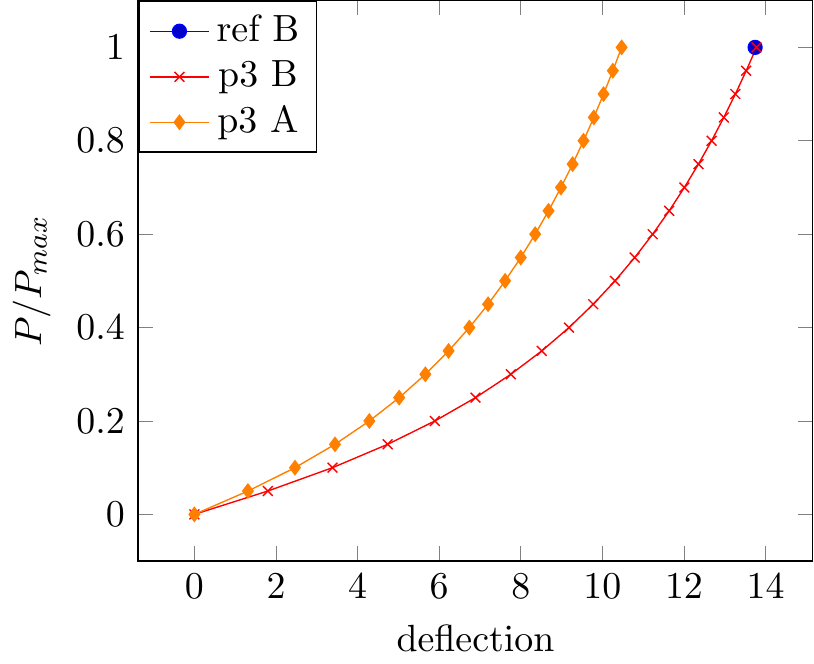}
	\caption{Vertical load-deflection for slit annular plate at points A and B with $10\times 80$ grid.}
	\label{fig:res_slit_ann2}
\end{figure}

The material and geometrical properties are $E=2.1\times10^8$, $\nu=0$, $R_i=6$, $R_o=10$, $t=0.03$ and $P_{max}=4.034$, see Figure \ref{fig:res_slit_ann_mesh}. We used structured quadrilateral meshes. The quantity of interest is the transverse displacement at point B. The reference value of $13.7432$ is taken from \cite{Ho01}. The initial and deformed mesh can be seen in Figure \ref{fig:res_slit_ann_def_mesh} and the results in Figure \ref{fig:res_slit_ann2} and Table \ref{tab:res_slitt_ann_end}.

\begin{table}[h!]
\begin{center}
\begin{tabular}{c|c|c}
 ref & p1 & p3\\
\hline
 13.7432 & 13.8224 & 13.7772\\	
\end{tabular}
\end{center}
\caption{Vertical deflection at point $B$ at maximal load for slit annular plate with $10\times 80$ grid.}
\label{tab:res_slitt_ann_end}
\end{table}

\subsection{Hemispherical shell subjected to alternating radial forces}
\label{subsec:hemi_shell}

\begin{figure}[h!]
	\centering
	\includegraphics[width=0.3\textwidth]{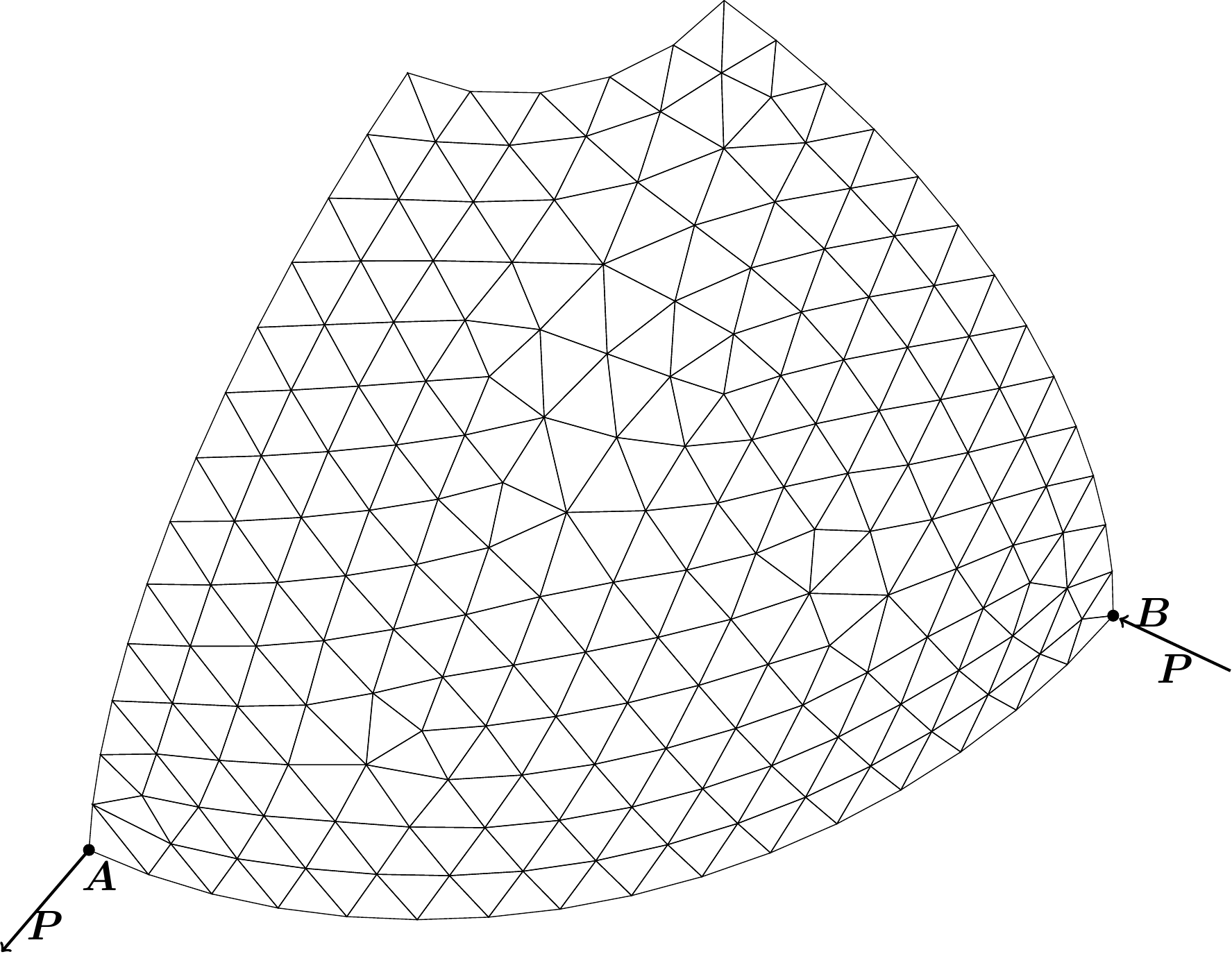}
	\caption{Geometry of hemispherical shell subjected to alternating radial forces with $h=1$.}
	\label{fig:hemi_shell_mesh}
\end{figure}

\begin{figure}[h!]
	\centering
	\includegraphics[width=0.24\textwidth]{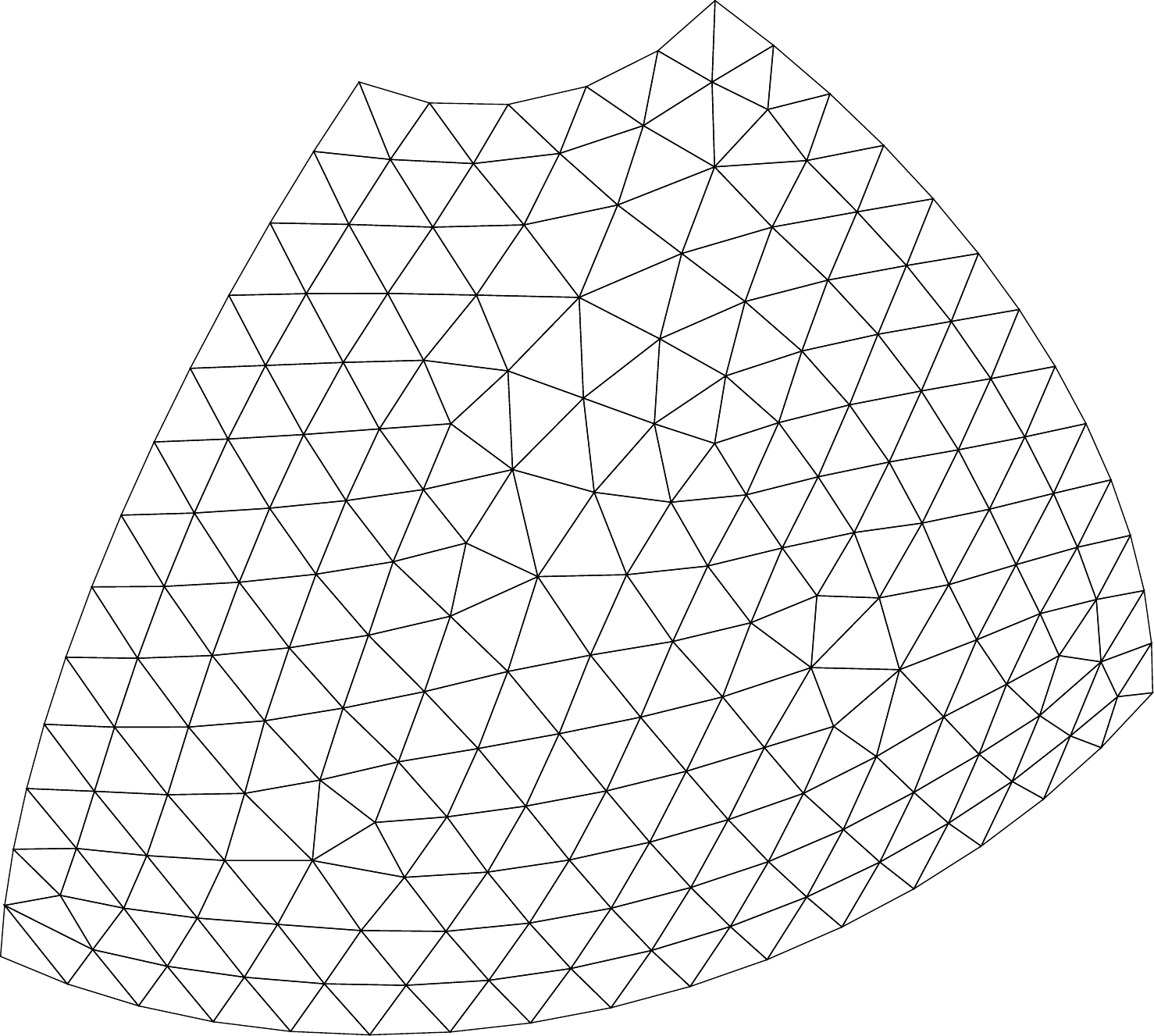}
	\includegraphics[width=0.24\textwidth]{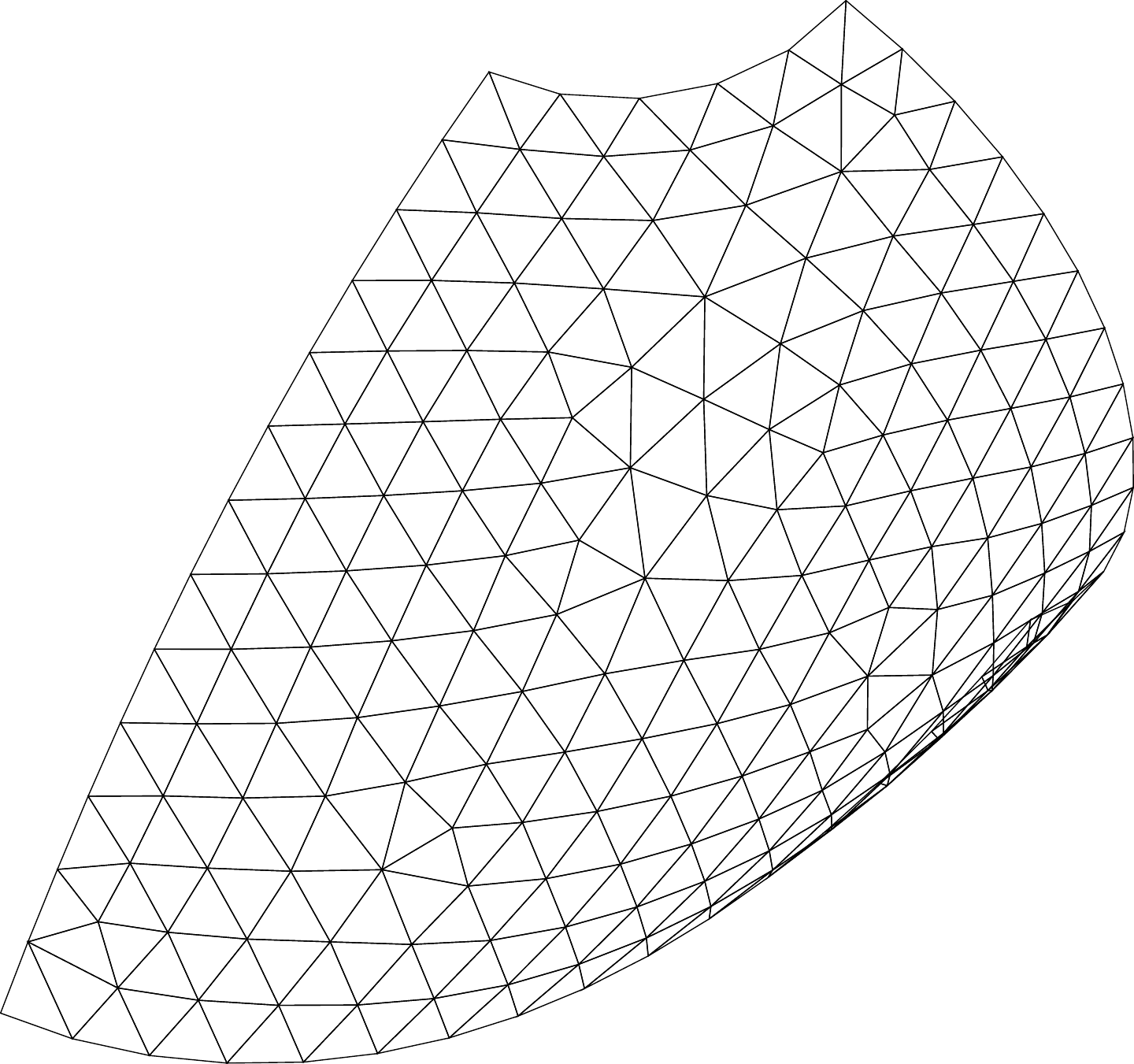}
	\caption{Initial and final configuration of hemispherical shell subjected to alternating radial forces with $h=1$.}
	\label{fig:hemi_shell_init_def_mesh}
\end{figure}

\begin{figure}[h!]
	\centering
	\includegraphics[width=0.45\textwidth]{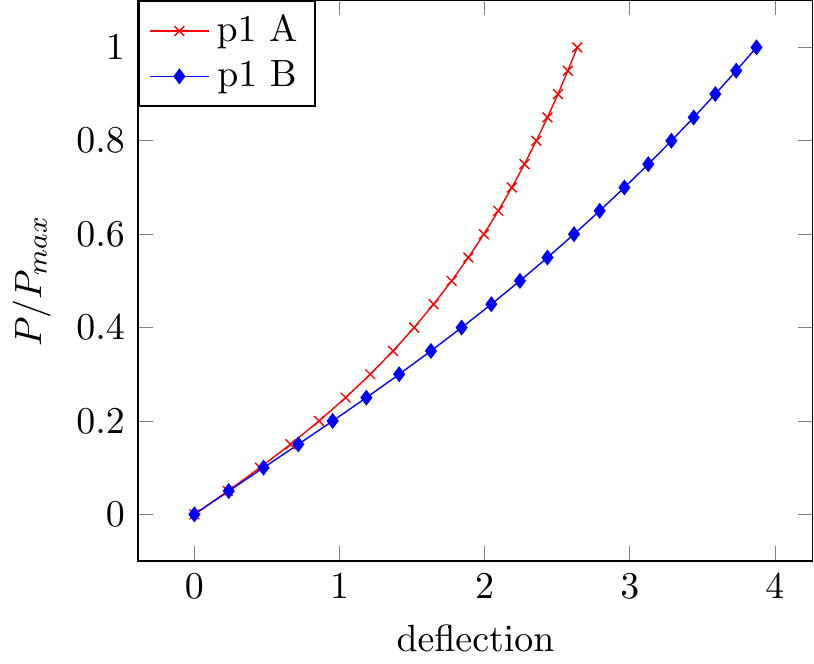}
	\caption{Radial load-deflections for the hemispherical shell subjected to alternating radial forces with mesh-size $h=0.25$.}
	\label{fig:res_hemisphere_alt}
\end{figure}

The material and geometrical properties are $E=6.825\times10^7$, $\nu=0.3$, $R=10$, $t=0.04$, see Figure \ref{fig:hemi_shell_mesh}. A non-structured triangulation is used with different mesh-sizes. For $P_{max}=1$ \cite{SIMO89} gives the reference value of the vertical deflection at point $B$ with $0.093$ at maximal load. In Table \ref{tab:hemi_shell_lin} the results for p1 and three different meshes can be found. For the large displacement case we used $P_{max}=400$, see Figure \ref{fig:hemi_shell_init_def_mesh} and \ref{fig:res_hemisphere_alt}. The results shown in Table \ref{tab:res_hemi_shell} are convenient with \cite{Ho01}.

\begin{table}[h!]
\begin{center}
\begin{tabular}{cccc}
ref & p1 1 & p1 0.5 & p1 0.25\\
\hline
0.093 & 0.110 & 0.092 & 0.0927
\end{tabular}
\caption{Radial load-deflection at point $B$ for the hemispherical shell subjected to alternating radial forces at maximal load for $P_{max}= 1$.}
\label{tab:hemi_shell_lin}
\end{center}
\end{table}
\begin{table}[h!]
\begin{center}
\begin{tabular}{c|cccc}
$h$& $2$ & $1$ & $0.5$ & $0.25$\\
\hline
p1 & 4.1218 & 3.8811 & 3.8560 & 3.8735\\
p3 & 3.8319 & 3.8781 & 3.8796 & 3.8796
\end{tabular}
\caption{Radial load-deflection at point $B$ for the hemispherical shell subjected to alternating radial forces at maximal load for $P_{max}=400$.}
\label{tab:res_hemi_shell}
\end{center}
\end{table}
\subsection{Twisted beam}
\label{subsec:twisted_beam}
\begin{figure}[h!]
	\centering
	\includegraphics[width=0.45\textwidth]{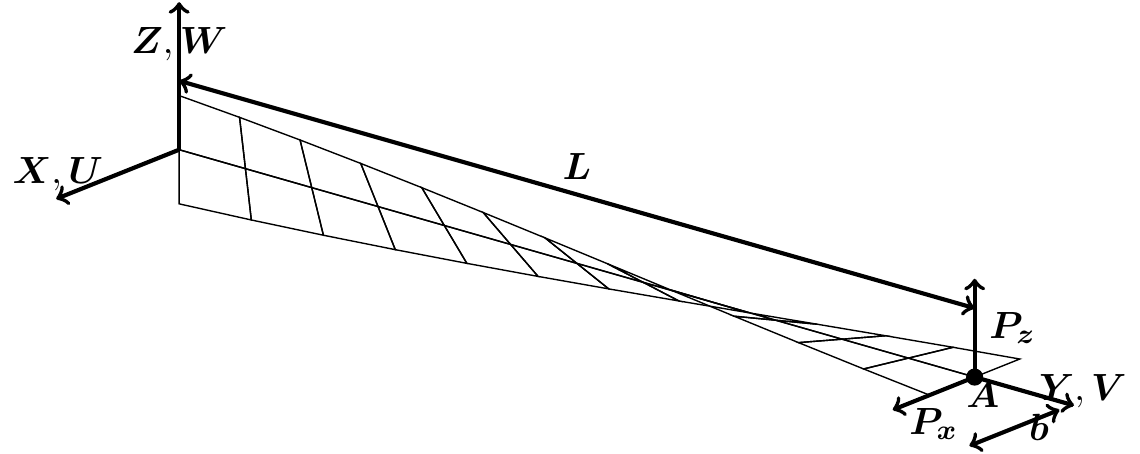}
	\caption{Geometry of twisted beam.}
	\label{fig:twisted_beam_mesh}
\end{figure}

A beam is twisted by $90$ degrees and clamped on the left side, whereas a point load is applied on the middle of the right boundary. The material and geometrical properties are $E=2.9\times10^7$, $\nu=0.22$, $L=12$, $b=1.1$, $t=0.0032$, $0.32$, see Figure \ref{fig:twisted_beam_mesh}. 

Different forces, $P_{max}\in\{10^{-6}$, $10^{-3}$, $1$, $10^3\}$, are applied in  x- and z-direction. Some combinations of thickness and force parameters led to a solution in a linear regime, see Table \ref{tab:res_twisted_beam0} and \ref{tab:res_twisted_beam2}, where the reference solutions are taken from \cite{BWS89} and \cite{MH85}, respectively. Others are already in the nonlinear regime, see Table \ref{tab:res_twisted_beam1} and \ref{tab:res_twisted_beam3}. Therefore, the full three-dimensional model is used with a $150\times14\times2$ structured cubic grid and standard Lagrangian elements of polynomial order 3, i.e. 162 dofs/cube, to generate a reference solution.

\begin{table}[h!]
	\begin{center}
		\begin{tabular}{ccc|cc}
			 & p1 & p3  & p1 & p3\\
			\hline
			2x12 & 5.504 & 5.460 & 1.332 & 1.339 \\
			4x24 & 5.470 & 5.452 & 1.337 & 1.338 \\
			6x36 & 5.460 & 5.451 & 1.337 & 1.338 \\
			8x48 & 5.456 & 5.451 & 1.337 & 1.338
		\end{tabular}
		\caption{Deflection $U_A\times10^3$ for $P_x =10^{-6}$, $P_z=0$, and $t=0.0032$ and  $W_A\times 10^3$ for $P_x =0$, $P_z=10^{-6}$, and $t=0.0032$ of twisted beam.  Reference values are $5.256$ and$1.294$. }
		\label{tab:res_twisted_beam0}
	\end{center}
\end{table}

\begin{table}[h!]
	\begin{center}
		\begin{tabular}{ccc|cc}
			& p1 & p3 & p1 & p3\\
			\hline
			2x12 & 4.555 & 4.536 & 1.337 & 1.282 \\
			4x24 & 4.538 & 4.532 & 1.279 & 1.280 \\
			6x36 & 4.535 & 4.532 & 1.280 & 1.280 \\
			8x48 & 4.534 & 4.531 & 1.280 & 1.280 
		\end{tabular}
		\caption{Deflection $U_A$ for $P_x =10^{-3}$, $P_z=0$, and $t=0.0032$  and $W_A$ for $P_x =0$, $P_z=10^{-3}$, and $t=0.0032$ of twisted beam. Reference values are $4.496$ and $1.227$.}
		\label{tab:res_twisted_beam1}
	\end{center}
\end{table}

\begin{table}[h!]
	\begin{center}
		\begin{tabular}{ccc|cc}
			 & p1 & p3 & p1 & p3\\
			\hline
			2x12 & 5.654 & 5.598 & 1.933 & 1.798 \\
			4x24 & 5.605 & 5.591 & 1.822 & 1.795 \\
			6x36 & 5.597 & 5.590 & 1.806 & 1.795 \\
			8x48 & 5.593 & 5.589 & 1.801 & 1.795
		\end{tabular}
		\caption{Deflection $U_A\times10^3$ for $P_x =1$, $P_z=0$, and $t=0.32$ and $W_A\times 10^3$ for $P_x =0$, $P_z=1$, and $t=0.32$ of twisted beam. Reference values are $5.424$ and $1.754$.}
		\label{tab:res_twisted_beam2}
	\end{center}
\end{table}

\begin{table}[h!]
	\begin{center}
		\begin{tabular}{ccc|cc}
			& p1 & p3  & p1 & p3\\
			\hline
			2x12 & 4.661 & 4.621 & 1.908 & 1.789 \\
			4x24 & 4.628 & 4.618 & 1.810 & 1.786 \\
			6x36 & 4.622 & 4.617 & 1.796 & 1.785 \\
			8x48 & 4.619 & 4.617 & 1.791 & 1.785
		\end{tabular}
		\caption{Deflection $U_A$ for $P_x =10^3$, $P_z=0$, and $t=0.32$ and $W_A$ for $P_x =10^3$, $P_z=1$, and $t=0.32$ of twisted beam. Reference values are $4.610$ and $1.778$.}
		\label{tab:res_twisted_beam3}
	\end{center}
\end{table}

\subsection{Z-section cantilever}
\label{subsec:zsection_cant}
\begin{figure}[h!]
	\centering
	\includegraphics[width=0.45\textwidth]{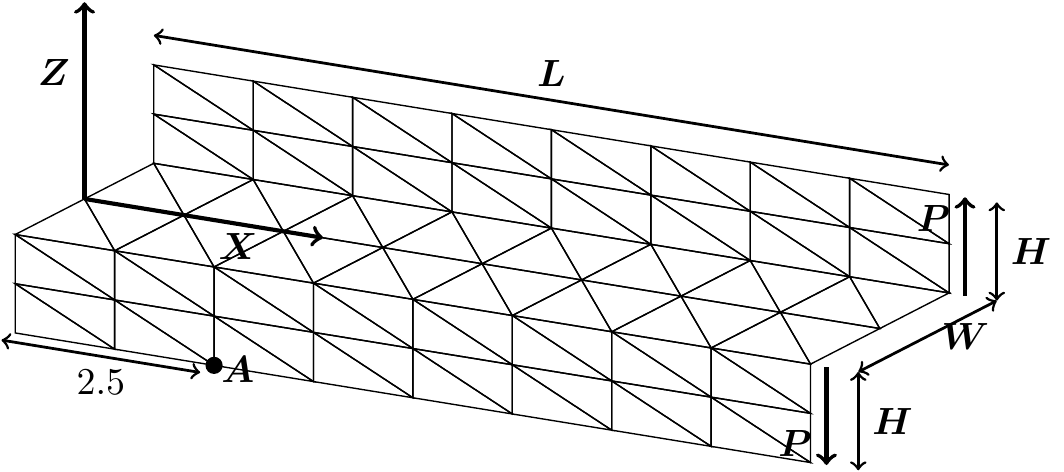}
	\caption{Geometry of Z-section cantilever.}
	\label{fig:zsec_cant_mesh}
\end{figure}

A moment $M=1.2\times 10^6$ is applied at the right end of a Z-section, which is fixed on the left side. Therefore, two shear forces $P=6\times 10^5$ are involved, see Figure \ref{fig:zsec_cant_mesh}. The material and geometrical properties are $E=2.1\times 10^{11}$, $\nu=0.3$, $t=0.1$, $L = 10$, $W=2$ and $H=1$. The quantity of interest is the membrane stress $\Sigma_{xx}$ at point $A$. The reference value $-1.08\times 10^8$ is taken from NAFEMS \cite{NAFEMS}. The results are compared with rotation-free elements \cite{FO07} and can be found in Table \ref{tab:res_zsec_cant}.

\begin{table}[h!]
\begin{center}
\begin{tabular}{cccc}
&\cite{FO07} & p1 & p3\\
\hline
8x6& $-0.953$ & $-0.7620$ & $-1.0929$ \\
32x15 & $-1.063$ & $-1.0777$  & $-1.0933$\\
64x30& - & $-1.0989$ & $-1.0933$\\
\end{tabular}
\caption{Membrane stress $\Sigma_{xx}\times 10^8$ of Z-section cantilever at maximal load.}
\label{tab:res_zsec_cant}
\end{center}
\end{table}

\subsection{T-section cantilever}
\label{subsec:tstruct}
\begin{figure}[h!]
\centering
\includegraphics[width=0.24\textwidth]{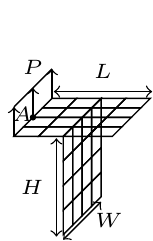}
\includegraphics[width=0.24\textwidth]{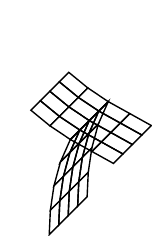}
\caption{Geometry of T-section cantilever and deformed configuration.}
\label{fig:Tsec_cant_mesh}
\end{figure}

We propose an example where more than two elements share an edge. The material and geometrical properties are $E=6\times10^{6}$, $\nu=0$, $t=0.1$, $L = 1$, $W=1$ and $H=1$. The structure is clamped on the bottom and a shear force $P_{max}=1000$ is applied on the left boundary, see Figure \ref{fig:Tsec_cant_mesh}.\newline

\begin{figure}[h!]
	\centering
	\includegraphics[width=0.45\textwidth]{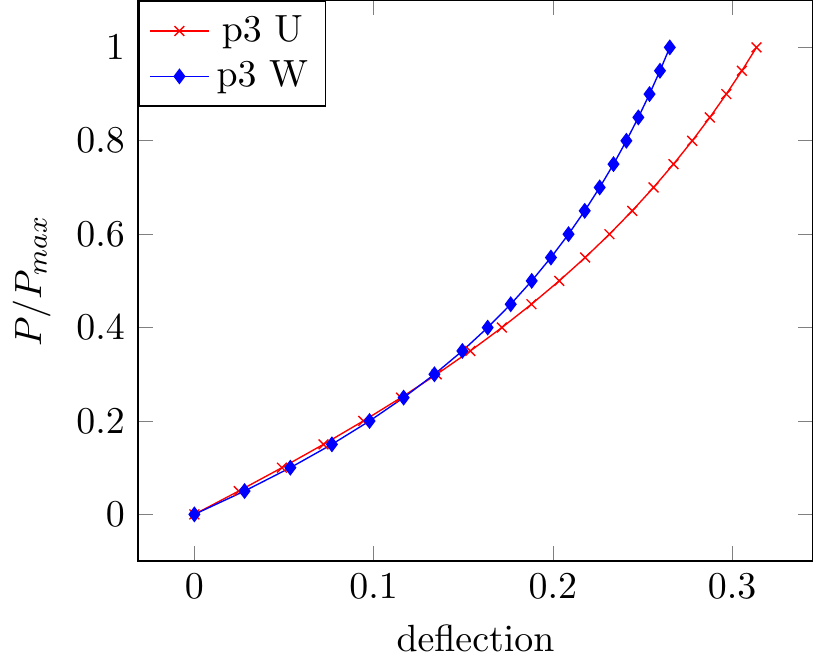}
	\caption{Horizontal and vertical deflection at point $A$ for T-section cantilever.}
	\label{fig:res_Tcant}
\end{figure}

The moment induced by the shear force $P$ on the left top branch goes over the kink to the bottom branch where the structure is fixed without inducing moments on the right top one. Thus, it only rotates and the curvature is zero also after the deformation. The deflections of the point $A$ are given in Figure \ref{fig:res_Tcant}.


\section*{Acknowledgments}
\label{sec:acknowledgments}
Michael Neunteufel has been funded by the Austrian Science Fund (FWF) project W1245.

\appendix
\section{Lagrange functional}
\label{sec:appendixA}
We compute the variations of the Lagrange functional in \eqref{eq:nonl_shell_method_first}, neglecting the sums over $\re{T}$ and $\re{E}$,
\begin{flalign}
\delta_{\physmom}\tilde{\Lagr} &= -\frac{12}{t^3}\langle\Mmat^{-1}\physmom,\delta\physmom\rangle+\langle\defgrad^T\refgrad(\physnv\circ\phi)-\refgrad\refnv,\delta\physmom\rangle_{\re{T}}&&\nonumber\\
&\quad-\langle\sphericalangle(\physnv_L, \physnv_R)\circ\phi-\sphericalangle(\refnv_L, \refnv_R),\delta\physmom_{\refbnv\refbnv}\rangle_{ \re{E}}\overset{!}{=}0,&&\label{eq:lag_var_mom}\\
\delta_{u}\tilde{\Lagr} &= \delta_u(\frac{t}{2}\|\Gstrain\|_{\Mmat})+\langle\physmom,\delta_u(\defgrad^T\refgrad(\physnv\circ\phi)-\refgrad\refnv)\rangle_{\re{T}}&&\nonumber\\
&\quad-\langle\delta_u(\sphericalangle(\physnv_L, \physnv_R)\circ\phi),\physmom_{\refbnv \refbnv}\rangle_{\re{E}}\overset{!}{=}0.&&\label{eq:lag_var_u}
\end{flalign}
Expressing $\physmom$ from \eqref{eq:lag_var_mom} and inserting it into \eqref{eq:lag_var_u} yields to the same expression as the variation of \eqref{eq:nonl_shell_method_first_u} with respect to the displacement $u$. We conclude that \eqref{eq:nonl_shell_method_first} and \eqref{eq:nonl_shell_method_first_u} are equivalent.\\

Equivalence of \eqref{eq:nonl_shell_method} and \eqref{eq:nonl_shell_method_first} follows by differentiating the identity $\defgrad^T\nu\circ\phi=0$ and some computations
\begin{flalign}
\langle\physmom,\defgrad^T\refgrad(\physnv\circ\phi)\rangle_{\re{T}} &= -\langle \mat{\defhesse_1:\physmom\\ \defhesse_2:\physmom\\ \defhesse_3:\physmom},\physnv\circ\phi\rangle_{\re{T}},&&
\end{flalign}
where
$\defhesse_i := \refgrad^2u_i+\refgrad((\Ptangref)_i)$, $(\Ptangref)_i$ denoting the i-th column of $\Ptangref$ and $\partial_{\re{x}_i}\defgrad$ the i-th partial derivative of $\defgrad$. With $\Ptangref = \Imat-\refnv\otimes\refnv$, neglecting $\phi$, and sum convention for $i$ we obtain
\begin{flalign}
\physnv\cdot\mat{\defhesse_1:\physmom\\ \defhesse_2:\physmom\\ \defhesse_3:\physmom}&= \physnv_i\refgrad((\Ptangref)_i+\refgrad^2 u_i):\physmom&&\nonumber\\
&=-\physnv_i(\refgrad(\refnv\otimes\refnv)_i-\refgrad^2u_i):\physmom&&\nonumber\\
&=-\physnv_i(\refgrad\refnv_i\otimes\refnv+\refnv_i\refgrad\refnv-\refgrad^2u_i):\physmom&&\nonumber\\
&=-(\physnv\cdot\refnv\refgrad\refnv-\physnv_i\refgrad^2u_i):\physmom&&\nonumber\\
&=-(\physnv\cdot\refnv\refgrad\refnv-\Hessian_{\physnv}):\physmom,&&\label{eq:split_defhesse}
\end{flalign}
where we used that $\refgrad\refnv_i\otimes\refnv:\physmom\equiv0$.
\section{Linearization}
\label{sec:appendixB}
To show that \eqref{eq:nonl_shell_method} simplifies to \eqref{eq:hhj_method} in the linear regime we use that the gradient of the displacement of the full three-dimensional model $\nabla U = \nabla(u+z\physnv\circ\phi)$ is small, $\nabla U =\mathcal{O}(\varepsilon)\ll 1$. Thus, we immediately obtain that $\refgrad u=\mathcal{O}(\varepsilon)$, $\defgrad=\Imat+\mathcal{O}(\varepsilon)$, $\Jbnd=1+\mathcal{O}(\varepsilon)$ and $\physmom=\mathcal{O}(\varepsilon)$. Furthermore, there holds $\physnv\circ\phi-\refnv = -\refnv^T\refgrad u + \mathcal{O}(\varepsilon^2)$ for $\varepsilon\rightarrow 0$. We neglect all terms of order $\mathcal{O}(\varepsilon^2)$ or higher.
For simplicity we will also neglect the $\phi$ dependency, e.g. we write $\physnv$ instead of $\physnv\circ\phi$.

Starting from \eqref{eq:final_B_plate}, we obtain on each $\re{T}\in\refT_h$
\begin{flalign}
\int_{\re{T}}\physmom:\Hessian_{\physnv} \,dx 
&= \int_{\re{T}}\sum_{i=1}^3\physmom:(\refgrad^2u_i\physnv_i) \,dx&&\nonumber\\
&\approx -\int_{\re{T}}\sum_{i=1}^3\physmom:\refgrad^2u_i(\refnv+\refnv^T\refgrad u)_i\,dx&&\nonumber\\
&\approx -\int_{\re{T}}\sum_{i=1}^3\physmom:\refgrad^2u_i\refnv_i\,dx&&\nonumber\\
&=\int_{\re{T}}\Div[\reftang]{\physmom}\cdot (\refnv^T\refgrad u)\,dx&&\nonumber\\
&\quad-\int_{\partial \re{T}} \refnv^T\refgrad u \cdot \physmom_{\refbnv}\,ds.&&\label{eq:lin_inner}
\end{flalign}
For the jump term we use \eqref{eq:approx_angle_comp} and \eqref{eq:approx_averaged_nv}, such that
\begin{flalign}
&\int_{\partial \re{T}}\frac{1}{2}\sphericalangle(\physnv_L, \physnv_R)\physmom_{\refbnv\refbnv}\,ds \approx\int_{\partial \re{T}}\overline{\Av{\physnv}}^n\cdot \physbnv\,\physmom_{\refbnv\refbnv}\,ds.&&\label{eq:jump_term_approx}
\end{flalign} 

For ease of presentation we neglect $\physmom_{\refbnv\refbnv}$ in \eqref{eq:jump_term_approx}, employ that $\frac{1}{\|\Ptauphys^{\perp}(\Av{\refnv})\|_2}\Ptauphys^{\perp}(\Av{\refnv})= \Ptauphys^{\perp}(\Av{\refnv})+\mathcal{O}(\varepsilon^2)$ and that $\Av{\refnv}=\refnv$ on a flat plane to obtain
\begin{flalign}
\int_{\partial \re{T}}\Ptauphys^{\perp}(\refnv)\cdot \physbnv\,ds &= \pm\int_{\partial \re{T}}\refnv\cdot(\physnv\times\physt)\,ds&&\nonumber \\
& \approx \mp\int_{\partial \re{T}}\refnv\cdot( \refnv^T\refgrad u\times\defgrad\reft)\,ds &&\nonumber\\
&\approx \mp\int_{\partial \re{T}}\refnv\cdot (\refnv^T\refgrad u\times\reft)\,ds&&\nonumber\\
&= \mp\int_{\partial \re{T}}\det(\refnv,\refnv^T\refgrad u,\reft)\,ds&&\nonumber\\
&= \mp\int_{\partial \re{T}}\det(\refnv,(\refnv^T\refgrad u\cdot \refbnv)\refbnv,\reft)\,ds&&\nonumber\\
&= \mp\int_{\partial \re{T}}\refnv^T\refgrad u\cdot \refbnv\underbrace{\det(\refnv,\refbnv,\reft)}_{=\mp 1}\,ds&&\nonumber\\
&= \int_{\partial \re{T}}\refnv^T\refgrad u\cdot \refbnv\,ds.&&\label{eq:lin_jump}
\end{flalign}

For $\int_{\partial\re{T}}\Av{\nu}\cdot\mu\,ds$ the linearization is done analogously and leads to the same result as \eqref{eq:lin_jump}.

If we now use \eqref{eq:lin_inner}, \eqref{eq:jump_term_approx}, \eqref{eq:lin_jump} and \eqref{eq:final_B_plate} and apply it to \eqref{eq:nonl_shell_method}, neglect the membrane energy term and the constants, and employ that $\refnv=\mat{0\\0\\1}$, we finally obtain
\begin{flalign}
-\|\physmom\|^2_{\Ltwo[\refshell_h]}+\sum_{\re{T}\in\refT_h}&\big(\int_{\re{T}}\refgrad w\cdot \Div[\reftang]{\physmom}\,dx&&\nonumber\\
&-\int_{\partial \re{T}}(\refgrad w)_{\reft} \physmom_{\refbnv\reft}\,ds\big),&&
\end{flalign}
which is indeed \eqref{eq:hhj_method}.\newline

\section{Variations}
\label{sec:appendixC}
We compute the variations of \eqref{eq:nonl_shell_method} to deduce the bilinear form of the according variational equations. Then we will (partly) integrate by parts to find the hidden boundary conditions in strong form.

For simplicity, we will neglect the material tensor $\Mmat$ and write only $\physnv$ instead of $\physnv\circ\phi$. The same holds for $\physbnv$ and $\physt$. We will consider only the formulation \eqref{eq:nonl_shell_method}, the case with the hybridization variable $\re{\alpha}$ in \eqref{eq:lagr_alpha} can be done analogously. 

Computing the first variation of  problem \eqref{eq:nonl_shell_method} with respect to $\physmom$ gives 
\begin{flalign}
&-\frac{12}{t^3}\langle\physmom,\delta\physmom\rangle_{}-\sum_{\re{T}\in\refT_h}\langle\delta\physmom,\Hessian_{\physnv} + (1-\refnv\cdot\physnv)\refgrad\refnv\rangle&&\nonumber\\
&-\sum_{\re{E}\in\refE_h}\int_{\re{E}}(\sphericalangle(\refnv_L, \refnv_R)-\sphericalangle(\physnv_L, \physnv_R) )\delta\physmom_{\refbnv\refbnv}\,ds=0&&\label{eq:var_mom}
\end{flalign}
for all permissible directions $\delta\physmom$. Testing \eqref{eq:var_mom} with functions which have only support on one edge $\re{E}$ of the triangulation $\re{\T}_h$ yields in strong
\begin{flalign}
\sphericalangle(\refbnv_L, \refbnv_R)-\sphericalangle(\physbnv_L, \physbnv_R)=0.&&\label{eq:angle_pres}
\end{flalign}

For the first variation of the membrane energy term of \eqref{eq:nonl_shell_method} in direction $v:=\delta u$ we immediately obtain for every $\re{T}\in\refT_h$
\begin{flalign}
\delta_u\|\Gstrain\|_{\re{T}}^2 &= \int_{\re{T}}(2\defgrad\Gstrain):\refgrad v\,dx.&&\label{eq_var_u_mem}
\end{flalign}
The other variations are more involved. We define the operator $\bar{(\cdot)}_{ij}:\R^{3\times3}\rightarrow\R^{2\times 2}$, which maps $3\times3$ matrices to its $2\times 2$ sub-matrix where the $i$-th row and $j$-th column are canceled out. Further, let  $\mathcal{A}^{ij}(\cdot):\R^{2\times 2}\rightarrow\R^{3\times 3}$ denotes the operator which embeds $2\times 2$ matrices into $3\times3$ matrices, such that $\overline{\mathcal{A}^{ij}(\Amat)}_{ij}=\Amat$ and the $i$-th row and the $j$-th column of $\mathcal{A}^{ij}(\Amat)$ are zero. Thus, $\mathcal{A}^{ij}(\cdot)$ is the right-inverse of $\bar{(\cdot)}_{ij}$.

With this, we define for $i,j\in\{1,2,3\}$
\begin{flalign}
\Amat^{\#}_{ij} := (-1)^{i+j}\mathcal{A}^{ij}\left(\cof{\bar{\Amat}_{ij}}\right).&&
\end{flalign}
Then, the following identity holds for all smooth matrix valued functions for $i,j,k\in\{1,2,3\}$
\begin{flalign}
\left(\frac{\partial}{\partial x_k}\cof{\Amat}\right)_{ij}= \Amat^{\#}_{ij}:\frac{\partial}{\partial x_k}\Amat.&&
\end{flalign}
With the notation $\Amat^{\#}_{\refnv,\physnv}:=\physnv_i\Amat^{\#}_{ij}\refnv_j$ and $v:=\delta u$ there further holds
\begin{flalign}
\delta_u J &= \defgrad^{\#}_{\refnv,\physnv}:\refgrad v,&&\\
\delta_u \physnv &= \frac{1}{J}\left(\mat{\defgrad^{\#}_{\refnv,1}:\refgrad v\\ \defgrad^{\#}_{\refnv,2}:\refgrad v\\ \defgrad^{\#}_{\refnv,3}:\refgrad v}-(\defgrad^{\#}_{\refnv,\physnv}:\refgrad v)\physnv\right),&&\label{eq:var_physnv}\\
\delta_u \Jbnd &=(\physt\otimes\reft):\refgrad v=\physt\cdot(\refgrad v)_{\reft}.&&
\end{flalign}
Note, that \eqref{eq:var_physnv} has the form of a covarient derivative.

By using $(\defgrad\reft)\times (\defgrad\refbnv)=\cof{\defgrad}\refnv$, $\defgrad^{\#}_{\refnv,i}:\refgrad v$ can be simplified to
\begin{flalign}
\defgrad^{\#}_{\refnv,i}:\refgrad v=\left((\refgrad v\, \reftang)\times(\defgrad\refbnv) + (\defgrad \reftang)\times (\refgrad v\,\refbnv)\right)_i.&&
\end{flalign}

Now, the volume term of \eqref{eq:final_B} is split into two terms depending on $u$
\begin{flalign}
\delta_u(\physmom:\Hessian_{\physnv})&= \physmom:\refgrad^2 v_i\physnv_i-\frac{1}{J}(\defgrad^{\#}_{\refnv,\physnv}:\refgrad v)\physmom:\Hessian_{\physnv}&&\nonumber\\
&\quad+\frac{1}{J}(\defgrad^{\#}_{\refnv,i}:\refgrad v)\,\physmom:\refgrad^2 u_i\label{eq:var_volterm1}&&
\end{flalign}
and
\begin{flalign}
\delta_u(-\refnv\cdot\physnv)=\frac{1}{J}\left((\refnv\cdot\physnv) \defgrad^{\#}_{\refnv,\physnv}-\defgrad^{\#}_{\refnv,\refnv}\right):\refgrad v.&&\label{eq:var_volterm2}
\end{flalign}
For the boundary integral of \eqref{eq:final_B} we use the averaged normal vector $\Av{\physnv}=\frac{1}{\|\physnv^o+\physnv\|_2}(\physnv^o+\physnv)$, with $\physnv^o$ denoting the element normal vector on the neighbored element. This yields
\begin{flalign}
&\delta_u(-\sphericalangle(\Av{\refnv}, \refbnv)+\sphericalangle(\Av{\physnv}, \physbnv))&&\nonumber\\
&=\frac{1}{\sqrt{1-(\Av{\physnv}\cdot\physbnv)^2}}\Big(\frac{\physnv^o\cdot\physnv}{\|\physnv+\physnv^o\|_2}(\delta_u\physnv)&&\nonumber\\
&\quad-\delta_u(\frac{1}{\|\physnv+\physnv^o\|_2}\physnv^o)\Big)\cdot\physbnv,&&\label{eq:var_bnd_term_dg}
\end{flalign}
which can be computed exploiting \eqref{eq:var_physnv}. Using \eqref{eq:approx_averaged_nv} instead of $\Av{\physnv}$ yields to a similar expression.

To obtain the boundary conditions of $u$ in strong form, which are hidden naturally in the weak form of the equation, we have to integrate by parts until no derivatives of $v$ appear.

E.g., \eqref{eq_var_u_mem} yields
\begin{flalign}
-\int_{\re{T}}2\Div[\reftang]{\defgrad\Gstrain}\cdot v\,dx+\int_{\partial\re{T}}2(\defgrad\Gstrain)_{\refbnv}\cdot v\,ds=0.&&
\end{flalign}

For \eqref{eq:var_volterm1} we have to integrate twice by parts obtaining
\begin{flalign}
&\int_{\re{T}}\Div[\reftang]{\Div[\reftang]{\physnv_i\physmom}}v_i+\Div[\reftang]{\frac{1}{J}(\physmom:\Hessian_{\physnv})\defgrad^{\#}_{\refnv,\physnv}}\cdot v&&\nonumber\\
&\qquad-\Div[\reftang]{\frac{1}{J}(\physmom:\nabla_{\reftang}^2u_i)\defgrad^{\#}_{\refnv,i}}\cdot v\,dx&&\nonumber\\
&\quad+\int_{\partial\re{T}}\Div[\reftang]{\physnv_i\physmom_{\refbnv}}v_i-\Div[\reftang]{\physnv_i\physmom}_{\refbnv}v_i &&\nonumber\\ &\qquad-\frac{1}{J}(\physmom:\Hessian_{\physnv}(\defgrad^{\#}_{\refnv,\physnv})_{\refbnv} - \physmom:\nabla_{\reftang}^2u_i(\defgrad^{\#}_{\refnv,i})_{\refbnv})\cdot v\,ds&&\nonumber\\
&\quad-\int_{\partial\partial\re{T}}\physnv\cdot v\,\physmom_{\refbnv\reft}\,dss=0,
\end{flalign}
where $\partial\partial T$ are the vertices of the element $T$ and $dss$ denotes point evaluation.

For \eqref{eq:var_volterm2} we get
\begin{flalign}
&-\int_{\re{T}}\Div[\reftang]{\frac{1}{J}\left((\refnv\cdot\physnv) \defgrad^{\#}_{\refnv,\physnv}-\defgrad^{\#}_{\refnv,\refnv}\right)}\cdot v\,dx&&\nonumber\\
&\quad+\int_{\partial\re{T}}\frac{1}{J}\left((\refnv\cdot\physnv) \defgrad^{\#}_{\refnv,\physnv}-\defgrad^{\#}_{\refnv,\refnv}\right)_{\refbnv}\cdot v\,ds=0.&&
\end{flalign}

Finally, one has to use integration by parts for \eqref{eq:var_bnd_term_dg} to obtain the last boundary terms. Adding up all boundary terms, taking care of the constants and material parameters, one obtain the natural boundary conditions in strong form with respect to the displacement $u$.

\bibliographystyle{acm}
\bibliography{zitate}
\end{document}